\documentclass[ijoo,nonblindrev]{informs-ijoo}

\OneAndAHalfSpacedXI

\usepackage[square,sort,comma,super,authoryear]{natbib}
\usepackage{etoolbox}
\usepackage{amssymb}
\TheoremsNumberedThrough     
\ECRepeatTheorems

\EquationsNumberedThrough    

\MANUSCRIPTNO{}

\usepackage{tikz}
\usetikzlibrary{hobby,backgrounds,calc,trees}
\usetikzlibrary{positioning}
\usepackage{cleveref}
\usepackage{subcaption}
\usepackage{bbm}
\begin{document}


 \RUNAUTHOR{Padmanabhan and Natarajan}

\RUNTITLE{Tree Bounds for Sums of Bernoulli Random Variables}

\TITLE{Tree Bounds for Sums of Bernoulli Random Variables: A Linear Optimization Approach}

\ARTICLEAUTHORS{%
\AUTHOR{Divya Padmanabhan}
\AFF{Engineering Systems Design, Singapore University of Technology and Design, \EMAIL{divya\_padmanabhan@sutd.edu.sg}} 
\AUTHOR{Karthik Natarajan}
\AFF{Engineering Systems Design, Singapore University of Technology and Design, \EMAIL{karthik\_natarajan@sutd.edu.sg}\footnote{The research was partly supported by the MOE Academic Research Fund Tier 2 grant T2MOE1706, “On the Interplay of Choice, Robustness and Optimization
in Transportation”.}}
} 

\ABSTRACT{%
We study the problem of computing the tightest upper and lower bounds on the probability that the sum of $n$ dependent Bernoulli random variables exceeds an integer $k$. Under knowledge of all pairs of bivariate distributions denoted by a complete graph, the bounds are NP-hard to compute. When the bivariate distributions are specified on a tree graph, we show that tight bounds are computable in polynomial time using linear optimization. These bounds provide robust probability estimates when the assumption of conditional independence in a tree structured graphical model is violated. Generalization of the result to finding probability bounds of order statistic for more general random variables and instances where the bounds provide the most significant improvements over univariate bounds is also discussed in the paper.
}%


\KEYWORDS{probability bounds, trees, linear optimization} \HISTORY{This paper was
first submitted in October 2019.}

\maketitle

%


\section{Introduction}
\label{sec:intro}
Analysis of the sums of Bernoulli random variables have received much attention among researchers in probability, computer science, optimization and engineering due to its wide applicability. For example, an insurer in risk management is interested in estimating the probability that the number of defaults among $n$ claims is $k$ or more \citep{Wang1998}. In the context of reliability, the probability that at least $k$ of $n$ subsystems is functional aids in estimating the probability that the entire system is functional \citep{boland1983}. In a retail environment, a popular measure of performance is the probability that there are stockouts at $k$ or more locations out of a total of $n$ locations given the inventory levels \citep{JiaXujie2012}.

Our interest is in the setting where bivariate dependence information among the Bernoulli random variables is available. Formally, denote by $\tilde{\bold{c}}$,  an $n$-dimensional Bernoulli random vector. Associated with this random vector is a graph $G = (V, E)$ where $V$ is the set of $n$ vertices denoting the random variables and $E$ is the set of edges denoting the pairs of random variables for which bivariate information is specified. The univariate distributions are denoted as $\mathbb{P}(\tilde{c}_i =r_i)$  for $i \in V$, $r_i \in \{0, 1\}$ and the bivariate distributions are denoted as $\mathbb{P}(\tilde{c}_i =r_i, \tilde{c}_j = r_j)$  for $(i,j) \in E$, $r_i \in \{0, 1\}, r_j \in \{0, 1\}$. Let $\Theta$ denote the set of distributions as follows:
\begin{align}
\Theta= \big\{
  \theta \in  \mathbb{P}(\{0,1\}^n) : & \ \mathbb{P}_\theta\left(\tilde{c}_i = 1, \tilde{c}_j = 1
    \right) = p_{ij} \text{ for }
    (i,j) \in E , \;\;  \mathbb{P}_\theta\left(\tilde{c}_i = 1
    \right) = p_i  \text{ for }
    i \in V  \big\},
\end{align}
where $\mathbb{P}(\{0,1\}^n)$ is the set of probability distributions of $n$ dimensional Bernoulli random vectors and $\mathbb{P}_\theta\left(\tilde{c}_i = 1, \tilde{c}_j = 0
    \right) = p_i - p_{ij}$, $\mathbb{P}_\theta\left(\tilde{c}_i = 0, \tilde{c}_j = 1
    \right) = p_j - p_{ij}$ and $\mathbb{P}_\theta\left(\tilde{c}_i = 0, \tilde{c}_j = 0\right) = 1- p_i - p_j + p_{ij}$.
Define  $ U(k)$ and $ L(k) $ as the largest and smallest possible probability that the sum of the $n$ random variables exceeds an integer $k$ computed over all distributions in $\Theta$:
\begin{align}
U(k)  &= \displaystyle \max_{\theta \in \Theta} \,\, \mathbb{P}_{\theta}\left(\sum_{i=1}^n \tilde{c}_i \geq k\right), \\
L(k)  &= \displaystyle \min_{\theta \in \Theta} \,\, \mathbb{P}_{\theta}\left(\sum_{i=1}^n \tilde{c}_i \geq k\right).
\end{align}
Unfortunately when the graph is complete with rational entries $p_{ij}$ and $p_i$, verifying if the set $\Theta$ is nonempty is an NP-complete problem (Theorem 3.3 in \citep{Pitowsky1991}). This implies that computing the tight bounds efficiently is highly unlikely, unless P = NP. For example, one can construct simple instances, even with $n = 3$ where the random variables are pairwise consistent (pairs of Bernoulli random variables exist) but not globally consistent (a multivariate Bernoulli random vector does not exist). An instance is $p_{1} = p_2 = p_3 = 1/2$ and $p_{12} = p_{23} = p_{13} = 0$ where $\Theta$ is empty while the bivariates are pairwise consistent (see \cite{Vorobev1962}).

\subsection{Tree Graphs}
A natural approach is to consider simpler graph structures where the feasibility of $\Theta$ is easy to verify. Towards, this, we consider the class of tree graphs $T = (V,E)$ which has attractive computational properties. Such graphs have been extensively studied in graphical models in computer science and machine learning (\citep{chowliutree1968,Lauritzen1996,Wainwright2008}) where a tree structured distribution that exploits conditional independence is used from the set $\Theta$.

We focus on a directed rooted tree $T$ representation of the graph where node $1$ is designated as the root and the arcs are directed away from the root node. Assume an arbitrary but fixed ordering of the remaining nodes.
The parent of a node $i \neq 1$ is denoted as $\text{par}(i)$ and is the unique node that connects to $i$ on the path from the root node. A child of a node $i$ is a node for which $i$ is the parent. A descendant of a node $i$  refers to any of the children of $i$ or the descendants of the children of $i$. A leaf node is a node with no descendants. We let $d_i$ denote the out-degree of node $i$ and denote the $s$th child of node $i$ (as per the ordering fixed a-priori) as $i(s)$. We denote by $T(i,s)$ the sub-tree rooted at $i$ consisting of the first $s$ sub-trees of $i$ where $V(i,s)$ is the set of vertices in $T(i,s)$ and $N(i,s)$ is the cardinality of this set. For ease of understanding, the notations are illustrated in Figure \ref{fig:tree-example} below.

\begin{figure}[h]
\caption{Suppose $n=7$ and the set of known bivariate marginal distributions corresponds to the set $E = \{(1,2),  (2,5), (2,6), (3,7), (3,1), (4,1) \}$. The figure gives the corresponding directed tree with all arcs pointing away from the root node $1$. The parent of nodes $5$ and $6$ is $2$, the parent of node $7$ is $3$ and the parent of nodes $2$, $3$ and $4$ is $1$. The degrees  of the various nodes are $d_1 =3, d_2 = 2, d_3 = 1, d_4 = 0, d_5 = 0, d_6 = 0, d_7 = 0$. Let $i$ denote the root node labelled 1. Assuming a non-decreasing order on the node labels, the  three children of $i$ are denoted as $i(1), i(2)$ and $i(3)$ (nodes 2,3 and 4 respectively). The sub-tree of $i$ induced by the first two children is denoted by $T(i,2)$ and is shaded. The number of nodes in $T(i,2)$ is denoted by $N(i,2)$. $N(i,2) = 6$ here. The set of vertices in $T(i,2)$ is $V(i,2) = \{ 1,2,3,5,6,7\}$.  }
\label{fig:tree-example}
\centering
\begin{tikzpicture}[edge from parent/.style={draw,-latex}]
 \tikzstyle{level 1}=[sibling distance=20mm]
\node (1) [label={right:$\leftarrow i$}] {1}
    child { node (2) [label={left:$i(1)\rightarrow$}]{ 2}
      child { node (5)  [label={below right:$ T(i,2)$}]  {5}
    }
      child { node (6) {6}
    }
  }
    child { node (3) [label={right: $\leftarrow i(2)$}]{3}
      child { node (7) {7}
    }
  }
   child { node (4) [label={right:$\leftarrow i(3)$}]{4}
   };

  \begin{pgfonlayer}{background}

  \draw[red,fill=gray,opacity=0.35](1.north west) to[closed,curve through={(5.north west).. (5.south west) .. (5.south east) ..(6.south) .. (7.south east).. (3.east) }] (1.north);
  \end{pgfonlayer}
\end{tikzpicture}
\end{figure}
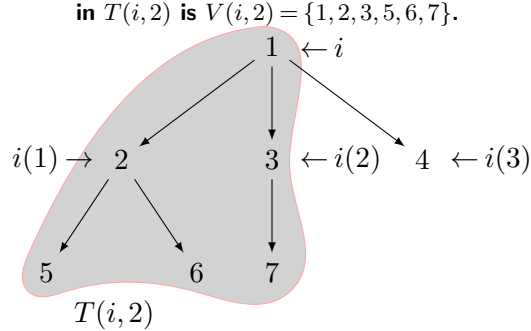
Given the bivariate distributions for this tree, a feasible distribution is given by:
\begin{align}
\label{eqn:cond-ind-whole-tree}
\mathbb{P}(\tilde{\bold{c}} = \bold{r}) = \mathbb{P}(\tilde{c}_1 = r_1) \prod_{j \neq 1}\mathbb{P}(\tilde{c}_j | \tilde{c}_{\text{par}(j)} = r_{\text{par}(j)}).
\end{align}
This distribution is based on the conditional independence assumption among the random variables in the tree which implies that for any two nodes $i \neq j$ such that $\text{par}(i) = \text{par}(j)$, we have:
\begin{align}
\mathbb{P}(\tilde{c}_i = x | \tilde{c}_{\text{par}(i)}, \tilde{c}_j) = \mathbb{P}(\tilde{c}_i =x | \tilde{c}_{\text{par}(i)}).\label{def:cond-ind}
\end{align}
As an example, the conditionally independent distribution for the tree shown in \Cref{fig:tree-example} is:
\begin{align}
\label{eqn:cond-ind-example}
\mathbb{P}(\tilde{\bold{c}} = \bold{r}) = & \mathbb{P}(\tilde{c}_1 = r_1)  \mathbb{P}(\tilde{c}_2 = r_2| \tilde{c}_1 = r_1 )  \mathbb{P}(\tilde{c}_3 = r_3| \tilde{c}_1 = r_1)  \mathbb{P}(\tilde{c}_4 = r_4| \tilde{c}_1 = r_1) \nonumber\\
& \,\,\, \times  \mathbb{P}(\tilde{c}_5 = r_5| \tilde{c}_2 = r_2)  \mathbb{P}(\tilde{c}_6 = r_6| \tilde{c}_2 = r_2)  \mathbb{P}(\tilde{c}_7 = r_7| \tilde{c}_3 = r_3)
\end{align}
Given a tree structured graphical model, many of the inference problems such as estimating the marginal distribution over a subset of random variables, computing the mode of the underlying distribution (see \citep{Lauritzen1996,Wainwright2008}) or estimating the probability that sum of the random variables exceed $k$ is easy. However, much lesser is discussed in the literature on what happens when the assumption of conditional independence is violated. In this paper, we use the tree structure of the graph as a certificate that the set $\Theta$ is nonempty and evaluate for the distributions in $\Theta$ that are extremal and attain the bounds $U(k)$ and $L(k)$. As our numerical results demonstrate, these bounds can in some cases be significantly different from the probability under conditionally independent distribution. Thus the bounds in the paper can provide robustness estimates on probabilities when the conditional independence assumption on the  underlying structured graphical models is violated. A similar problem was recently studied by \cite{dhara2019} where tight tree bounds were proposed for the expectation of a sum of discrete random variables beyond a threshold using linear optimization. In contrast to their work, our focus in this paper is on probability bounds which requires the use of different proof techniques.

\subsection{Related Results}
When only univariate probabilities $p_i$ are known for Bernoulli random variables, the tight upper bound on $\mathbb{P}\left(\sum_{i=1}^n \tilde{c}_i \geq k\right)$ was derived by \citep{ruger1978,morgenstern1980} as follows:
\begin{align}
\label{eqn:univar-bound}
\displaystyle U_{uv}(k) = \min \left( \left(\min_{t \in \{0, \ldots, k-1\}} \sum_{i=1}^{n-t} \frac{p_{(i)}}{k-t} \right), 1 \right)
\end{align}
where $p_{(1)} \leq p_{(2)} \leq \ldots \leq p_{(n)}$ are the order statistics of the probabilities $p_i$. When instead of exact bivariate probabilities, only lower bounds on the bivariate probabilities are known, a tight upper bound on the tail probability is computed in polynomial time for $k=1$ (union of events) in \cite{boros2014}. Verifying the existence of a feasible distribution in this case is possible in polynomial time (see Chapter 8 of \cite{bertsimas1997introduction}). When the exact probabilities $p_{ij}$ on all edges of a graph are known, the tight bound is obtained as a solution to an exponential sized linear programming formulation as discussed in \citep{hailperin1965,kounias1976,prekopa1997}. \cite{hunter1976} and \cite{worsley1982} proposed an upper bound for $\mathbb{P}(\sum_{i=1}^n c_i \geq 1)$ in terms of the total weight of a maximum spanning tree on a complete graph of $n$ vertices, with the weight of edge $(i,j)$ taken as the probability $p_{ij}$. Their proposed upper bound is $\sum_{i=1}^n p_i - \max_T \sum_{(i,j) \in T} p_{ij}$, where the maximum is computed over  all possible trees $T$. In the specific case where the bivariate probabilities are given as $p_{ij} = 0$ for all edges $(i,j)$ not in a tree $T$,  \cite{kounias1968} show that $\sum_{i=1}^n p_i - \sum_{(i,j) \in T} p_{ij}$ is a tight upper bound for $k = 1$.  Extensions of the approach to higher order information have been considered in \cite{tomescu1986hypertrees, bukszar2002hypercherry}. 
For tree structured bivariate information, \cite{ruschendorf1991} proposed a conditioning method for series and star graphs. \cite{Embrechts2010} also proposed upper bounds building on these results. These bounds are tight in very special cases and are in general not tight. 
We provide a snapshot of the results in \Cref{tab:positioning_of_work}.
\begin{table}
\centering
\caption{Tight upper bound on $\max \mathbb{P}(\sum_{i=1}^n \tilde{c}_i \geq k)$ for Bernoulli random variables}
\label{tab:positioning_of_work}
\begin{tabular}{|p{2cm}|p{4.5cm}|p{5cm}|p{2.5cm}|}
\hline
\textbf{Univariate} &\textbf{Bivariate} & \textbf{Solution approach} & \textbf{Computation} \\\hline
$p_i$  &  Not given &  Closed form bound \citep{ruger1978,morgenstern1980} & Easy \\\hline
$p_i$  & Lower bounds on bivariate probabilities in a complete graph &  Tight bound for $k = 1$ \citep{boros2014} & Easy \\\hline
$p_i$   &Exact values for bivariate probabilities in a complete graph &  Exponential sized linear program \citep{hailperin1965} & Hard\\\hline
$p_i$   & Exact values for bivariate probabilities in a complete graph; bivariate probabilities are $0$ for edges not in a tree &   Tight bound for $k=1$ \citep{hunter1976,worsley1982} & Easy  \\\hline
$p_i$ & Exact values for bivariate probabilities in a tree &  Linear program [Current paper]& Easy \\\hline
\end{tabular}
\end{table}

\subsection{Overview of Approach}
We consider the exponential sized linear program to compute $U(k)$ for a given graph:
\begin{align}
\displaystyle \max\limits_{\theta} & \sum\limits_{\bold{c} \in \{0,1\}^n} \theta(\bold{c})   \;\;\mathbbm{1}\lbrace \sum_i c_i \geq k \rbrace  \nonumber  \\
\textrm{s.t.} & \sum\limits_{\bold{c}: c_i = 1}  \theta(\bold{c})= p_i \; \text{ for } i \in V\nonumber \\
 & \sum\limits_{\bold{c}: c_i = 1, c_j = 1 }   \theta(\bold{c})= p_{ij} \; \text{ for } (i,j) \in E \nonumber \\
 & \sum\limits_{\bold{c} \in \{0,1\}^n} \theta(\bold{c})= 1 \nonumber \\
 &  \theta(\bold{c}) \geq 0 \;\;\; \text{ for } \bold{c} \in \{0,1\}^n \nonumber,
\end{align}
where $\mathbbm{1}\lbrace \sum_i c_i \geq k \rbrace = 1$ if $c_i \geq k$ and $0$ otherwise and $\theta(\bold{c})$ denotes the probability of realization $\bold{c}$. The first two constraints enforce the given information on the univariate and bivariate probabilities while the last two constraints ensure that $\theta$ is a valid distribution. The formulation above is exponential sized owing to number of realizations of $\tilde{\bold{c}}$.

The dual to the above formulation  is:
\begin{align}
\min\limits_{\lambda, \boldsymbol{\alpha},  \boldsymbol{\beta}}  &\ \; \lambda +  \sum_{i=1}^n \alpha_i p_i + \sum\limits_{(i,j) \in E} \beta_{ij} p_{ij} \nonumber \\
\textrm{s.t.} & \;
 \lambda + \sum_{i=1}^n \alpha_i c_i + \sum\limits_{(i,j) \in E} \beta_{ij} c_i c_j \geq 1,   \;\; \text{ for } \bold{c} \in \{0,1\}^n \text{ where} \sum\limits_{i=1}^n c_i \geq k \label{constr:1} \\
& \; \lambda + \sum_{i=1}^n \alpha_i c_i + \sum\limits_{(i,j) \in E} \beta_{ij} c_i c_j \geq 0 \;\; \text{ for }  \bold{c} \in \{0,1\}^n  \label{constr:0}
\end{align}
The dual has an exponential number of constraints and these can be  grouped into two sets of constraints. In particular, for the separation version of the above problem, given $\lambda, \boldsymbol{\alpha}, \boldsymbol{\beta}$, we need to verify if all constraints in \eqref{constr:1} and \eqref{constr:0} are met, or else we need to find a violated inequality. By the equivalence of separation and optimization in  \cite{grotschel1988geometric}, a polynomial time solution to the separation problem would imply  a polynomial time algorithm for the optimization problem. As we discuss in this paper, the separation problem for the dual is efficiently solvable when the graph is a tree and in particular, we can develop a compact linear program to compute $U(k)$.

In \Cref{sec:ccqkp} we consider the special case of a quadratic knapsack problem with cardinality constraints on a tree graph which arises in the dual formulation and develop a compact linear program based on a set of dynamic programming recursions. Building on this, in \Cref{sec:prob_bounds_tree}  we propose a polynomial sized linear programming formulation to compute $U(k)$ for tree graphs. We also compare this formulation with the corresponding approach for the conditionally independent distribution in a graphical tree model. We discuss generalization of the bounds to weighted sums of probabilities and orders statistics in Section 4 and provide numerical results in Section 5.

\section{Special Case of Cardinality Constrained Quadratic Knapsack}
\label{sec:ccqkp}
Given a graph $G = (V,E)$ and the parameter vectors $\boldsymbol{\alpha}, \boldsymbol{\beta}$, consider the quadratic optimization problem:
\begin{align}
\label{opt:card-constr-quad-knapsack}
\min \left\lbrace \sum_{i=1}^n \alpha_i c_i + \sum_{(i,j) \in E} \beta_{ij} c_i c_j:  \sum_{i=1}^n c_i \geq k, c_i \in \{0,1\} \forall \; i \in V \right\rbrace,
\end{align}
which arises in the dual formulation. Formulation (\ref{opt:card-constr-quad-knapsack}) is a special case of a quadratic knapsack problem with only cardinality constraints. The cost minimization version of a quadratic knapsack problem in its general form requires to find a vector $\bold{c} \in \{0, 1\}^n$ to minimize $  \sum_{i \in V} \alpha_i c_i + \sum_{(i,j) \in E} \beta_{ij} c_i c_j$ subject to  a constraint $\sum_{i \in V} a_i c_i \geq k$. Each vertex can be interpreted as corresponding to an item and each item is associated with a utility $a_i$. The requirement is to choose a set of items so that the overall utility of the selection is at least $k$. To this end if an item $i$ is selected, $c_i = 1$ else $c_i =0 $. Given an edge $(i,j)$ in the graph, an additional cost of $\beta_{ij} $ is incurred when both items $i$ and $j$ are selected, in addition to individual item costs $\alpha_i$ and $\alpha_j$. When item $i$ is selected but item $j$ is not selected, a cost of $\alpha_i$ alone is incurred for the item. The overall goal is therefore to select a set with a total utility of at least $k$ while minimizing the cost induced by the selected items.


Quadratic knapsack problem is  NP-hard in the strong sense \citep{Caprara1999,Fomeni2014}. However for special types of graphs such as series-parallel graphs (of which trees are a special case), a pseudo-polynomial time dynamic programming algorithm of time complexity $O(nk^2)$ is available \citep{RaderJr2002}.
 Our interest is on a special instance of the quadratic knapsack, where the graph $G$ is a tree, $a_i = 1$ for all $i$ and at least $k$ items need to be selected. For this problem, when the graph is a tree, a dynamic programming algorithm has been proposed in \cite{Billionet-cardQP}.  In this section, we will develop a linear optimization formulation for this problem that builds on dynamic programming. The advantage of the linear optimization formulation is that it can in turn be used in the computation of the probability bounds where the parameters $\boldsymbol{\alpha}, \boldsymbol{\beta}$ are themselves decision variables to develop a compact linear program.

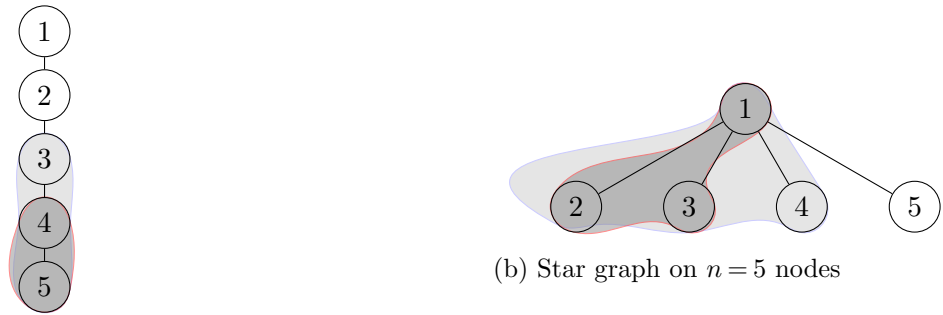
\begin{figure}[h!]
\begin{subfigure}{0.5\textwidth}
\centering
\begin{tikzpicture}[
  every node/.style = {shape=circle,
    draw, align=center},level distance=0.85cm]]
  \node {1}
    child  { node (2) {2}
    child  { node (3) {3}
      child  { node (4) {4}
       child  { node (5) {5}} }}};
    \begin{pgfonlayer}{background}
  \draw[blue,fill=gray,opacity=0.2](3.east) to[closed,curve through={ (4.east).. (5.east).. (5.south east) .. (5.south).. (5.south west) .. (4.west) .. (3.west)}] (3.north);
      \draw[red,fill=gray,opacity=0.45](4.east) to[closed,curve through={(5.east).. (5.south east) .. (5.south west) .. (4.west) }] (4.north);
  \end{pgfonlayer}

\end{tikzpicture}
\caption{Series graph on $n=5$ nodes}
\label{fig:series-graph-5}
\end{subfigure}
\begin{subfigure}{0.5\textwidth}
\centering
\begin{tikzpicture}[
  every node/.style = {shape=circle,
    draw, align=center,
 },level distance=1.3cm]]
  \node (1) {1}
    child { node (2) {2} }
    child { node (3) {3}}
    child { node (4) {4}}
      child { node (5) {5} };

        \begin{pgfonlayer}{background}
          \draw[blue,fill=gray,opacity=0.2](1.north) to[closed,curve through={(1.east)..(4.north east).. (4.east).. (4.south) .. (4.south west).. (3.south) .. (2.south east) .. (2.south west).. (1.west) }] (1.north west);
  \draw[red,fill=gray,opacity=0.45](1.north) to[closed,curve through={(1.east)..(3.north east).. (3.east).. (3.south) .. (3.south west).. (2.south).. (2.west) .. (2.north west) .. (1.west) }] (1.north west);
  \end{pgfonlayer}
\end{tikzpicture}
\caption{Star graph on $n=5$ nodes}
\label{fig:star-graph-5}
\end{subfigure}
\caption{Illustrative examples representing the basic components of any general tree. The shaded parts illustrate the relevant sub-trees corresponding to smaller sub-problems.}
\label{fig:series_star_trees}
\end{figure}

We start with the linear programming formulations for two particular trees - the series graph and the star graph (see Figure \ref{fig:series_star_trees}) and generalize the result to arbitrary trees by viewing them as a combination of several series and/or star sub-graphs. Throughout the paper, we use the notation $[n]$ to denote the set $\{1, \ldots, n \}$ for any integer $n$, and  $[i,j]$ to denote the set $\{ i, i+1, \ldots, j\}$ for integers $i$ and $j$.

\subsection{Formulation for a Series Graph}

The series graph $G_{series}=(V, E_{series})$ on $n$ nodes is a graph with the edge set $E_{series} =\{ (1,2), (2,3), \ldots, (n-1,n) \}$.  Such a graph contains exactly one leaf node (node $n$) (see  \Cref{fig:series-graph-5}). We will now propose a linear programming  formulation for solving the quadratic minimization problem with a series graph:
\begin{align}
\label{opt:card-constr-quad-knapsack-series}
\bar{Q} = \min \left\lbrace \sum_{i=1}^n \alpha_i c_i + \sum_{i \in [n-1]} \beta_{i,i+1} c_i c_{i+1}:  \sum_{i=1}^n c_i \geq k, c_i \in \{0,1\} \forall \; i \right\rbrace
\end{align}

\begin{proposition}
\label{prop:series_graph_lp_ccqkp}
For a series graph $G_{series} = (V, E_{series})$,  the optimal value $\bar{Q}$ of the quadratic knapsack problem in \eqref{opt:card-constr-quad-knapsack-series} can be obtained by solving the following linear program:
\begin{align*}
\begin{array}{rllll}
\bar{Q} =   \max\limits_{z, \bold{x}} & z & \\
 \mbox{s.t.} & z \leq x_{1, 0, t},  & t \in [k, n-1]  \\
&  z  \leq x_{1, 1, t}  & t \in [k, n]  \\
&  x_{i,0,t}  \leq x_{i+1,0,t}, &  t \in [0,n-i-1], i \in [1,n-1] \\
&   x_{i,0,t} \leq x_{i+1,1,t}, & t \in [1,n-i] , i \in [1,n-1]  \\
&  x_{i,1,t}\leq x_{i+1,0,t-1} + \alpha_i, & t \in [1, n-i] , i \in [1,n-1]   \\
&  x_{i,1,t}  \leq x_{i+1,1,t-1} + \alpha_i  + \beta_{i,i+1}, & t \in [2, n-i+1], i \in [1,n-1] \\
&  x_{n,0,0}  = 0 & \\
&  x_{n,1,1} = \alpha_n. &
\end{array}
\end{align*}
\end{proposition}
\noindent \textit{Proof:} Denote by $\bar{f}(i,y,t)$ the optimal value of formulation (\ref{opt:card-constr-quad-knapsack-series}), when restricted to the sub-tree rooted at $i$ such that $t$ nodes are selected from this sub-tree and $c_i$ takes a value $y \in \{0,1\}$:
\begin{align}
\label{opt:dp-defn-f-series}
\bar{f}(i,y,t) = \min \left\lbrace \sum\limits_{l \geq i } \alpha_l c_l + \sum\limits_{l \geq i }  \beta_{l,l+1} c_l c_{l+1}:  \sum_{l \geq i} c_l = t, c_i = y, c_l \in \{0,1\} \forall l \geq i  \right\rbrace
\end{align}
If $y = 0$, the range of admissible values for $t$ is $[0,n-i]$ while if $y=1$, the range of admissible values for $t$ is $[1, n-i+1]$ (see the shaded region in \Cref{fig:series-graph-5} for the relevant sub-tree for $\bar{f}(3, \cdot, \cdot)$).

The smallest such sub-tree is rooted at node $n$ which contains the leaf node $n$ alone. With $c_n = y$, the only possible value of $t$ is $y$ which leads to the following two base-cases:
\begin{align}
\label{rec:series-base}
 \bar{f}(n,0,0) = 0,   \bar{f}(n,1,1) = \alpha_n.
 \end{align}
Using the above two base cases, the optimal value for other sub-problems can be recursively computed. We will now develop the recursions for $\bar{f}(i,y,t) $ for any internal node $i$, for all valid values of $y$ and $t$ in terms of the optimal values for $\bar{f}(i+1,\cdot, \cdot) $.
For any valid value of $t$, we have:
\begin{align}
\bar{f}(i,0,t) &= \min( \bar{f}(i+1,0,t), \bar{f}(i+1,1,t)) \,\,  \label{rec:series1} \\
\bar{f}(i,1,t) &= \min( \bar{f}(i+1,0,t-1) + \alpha_i,\bar{f}(i+1,1,t-1) + \alpha_i  + \beta_{i,i+1}) \,\,   \label{rec:series3}
\end{align}
When $c_i = 0$, the $t$ nodes which take a value of $1$ must all be located in the sub-tree rooted at $i+1$.
The terms in the right hand side of \Cref{rec:series1} deal with the case where $c_{i+1} = 0 $ and $c_{i+1} = 1$ respectively. On the other hand, when $c_i = 1$, then the sub-tree rooted at $i+1$ must select $t-1$ nodes so that a total of $t$ nodes are selected from the subtree rooted at $i$. Here if $c_{i+1}$ takes a value of 0, only an additional cost of $\alpha_i$ is incurred, while if $c_{i+1} = 1$, then an additional cost of  $\alpha_i + \beta_{i,i+1}$ is incurred.   In all these cases the range of valid $t$ varies depending on the values of $c_i$ and $c_{i+1}$.  For example, if $c_i= 0$ and $c_{i+1} = 1$, the range of permissible values of $t$ is $[1, n-i]$ while for $c_i= 1$ and $c_{i+1} = 1$, the permissible range is $[2, n-i+1]$.

The optimal objective $\bar{Q}$ on the overall series graph is obtained by looking at the optimal values of $\bar{f}(\cdot, \cdot, \cdot)$ corresponding to the root node. In particular,
\begin{align}
 \bar{Q} &= \displaystyle \min\left(\min_{t_1 \in [k, n-1]} \bar{f}(1, 0, t_1),\min_{t_2 \in [k,n]}\bar{f}(1, 1, t_2) \right).  \label{rec:root1}
\end{align}
While the range of permissible values of $t_1$ in $\bar{f}(1, 0, t_1)$ is  $[k, n-1]$ (as $c_1  = 0$ here), the permissible range of $t_2$ is $[k,n]$ in  $\bar{f}(1, 1, t_2)$.
The variable $\bold{x}$ which is $O(n^2)$ encodes the optimal values $\bar{f}(\cdot,\cdot,\cdot)$ of all the sub problems while $z$ encodes the optimal value $\bar{Q}$. The inequalities in the linear program arise as a consequence of linearizing the minimum operator in \Cref{rec:root1,rec:series1,rec:series3}.
\halmos

\subsection{Formulation for a Star Graph}
The star graph $G_{star}=(V, E_{star})$ on $n \geq 2$ nodes is a graph with edge set $E_{star} =\{ (1,2), (1,3), \ldots, (1,n) \}$ (see \Cref{fig:star-graph-5} for an illustration). We will now propose a linear programming formulation for solving the minimization problem:
\begin{align}
\label{opt:card-constr-quad-knapsack-star}
Q_{star} = \min \left\lbrace \sum_{i=1}^n \alpha_i c_i + \sum_{i \in [2,n]} \beta_{1i} c_1 c_{i}:  \sum_{i=1}^n c_i \geq k, c_i \in \{0,1\} \forall i \right\rbrace
\end{align}

\begin{proposition}
\label{prop:star_graph_lp_ccgkp}
For a star graph $G_{star} = (V, E_{star})$, the optimal value $Q_{star}$ of the quadratic knapsack problem in \eqref{opt:card-constr-quad-knapsack-star} can be obtained by solving the following linear program:
\begin{align*}
\begin{array}{rllll}
Q_{star} = \max\limits_{z, \bold{x}} & z \\
\mbox{s.t.} & z \leq x_{n, 0, t}, & t \in [k, n-1]  \\
 & z  \leq x_{n, 1, t}, & t \in [k, n]  \\
& x_{i,0,t}  \leq x_{i-1,0,t}, & t \in [0,i-2], i \in [3,n] \\
& x_{i,0,t} \leq x_{i-1,0,t-1} + \alpha_i, &  t \in [ 1, i-1] , i \in [3,n]   \\
& x_{i, 1, t}  \leq x_{i-1,1,t}, &  t \in [1, i-1], i \in [3,n]  \\
& x_{i,1,t} \leq x_{i-1,1,t-1} + \alpha_i + \beta_{1i}, & t \in [2, i] , i \in [3,n]\\
& x_{2,0,0} = 0,\,\,\, x_{2, 0, 1} = \alpha_2 & \\
& x_{2, 1, 1} = \alpha_1, \,\,\, x_{2, 1, 2} = \alpha_1 + \alpha_2 + \beta_{12}.
\end{array}
\end{align*}
\end{proposition}
\textit{Proof:}
For values of $i \geq 2$, denote by $g(i,y, t)$ the minimum value that is obtained by restricting attention to the sub-tree containing the nodes $\{ 1, \ldots, i \}$, such that $t$ nodes are selected from this sub-tree and $c_1$ takes a value $y$:
\begin{align}
\label{opt:dp-defn-g-star}
g(i,y, t) = \min \left\lbrace \sum\limits_{l =1 }^i \alpha_l c_l +  \sum_{l=2}^i \beta_{1, l} c_1 c_{l}:  \sum_{l = 1}^i c_l = t, c_1 = y, c_l \in \{0,1\} \forall l \leq i  \right\rbrace
\end{align}
The region shaded in \Cref{fig:star-graph-5} (consisting of nodes 1,2,3,4) shows the relevant tree for $g(4, \cdot, \cdot)$. If $c_1= 0$, the valid values of $t$ lie in $[0, i-1]$ while if $c_1 = 1$, $t$ must lie in $[1, i]$.

We will now provide recursions to compute the values of $g(\cdot,\cdot,\cdot)$. The base conditions look at the sub-tree with exactly two nodes $\{1, 2 \}$ as follows:
\begin{align}
g(2,0,0) = 0, g(2, 0, 1) = \alpha_2 \label{star-base-cond-c1-0}\\
g(2, 1, 1) = \alpha_1, g(2, 1, 2) = \alpha_1 + \alpha_2 + \beta_{12}\label{star-base-cond-c1-1}
\end{align}
\Cref{star-base-cond-c1-0} deals with the case  where $c_1 = 0$ and the possible value of $t$ is either $0$ or $1$,  depending on the value of $c_2$. If $t=0$,  it must be that $c_2 = 0$ in which case, no cost is incurred as none of the nodes are selected. If $t=1$, then the only possibility is $c_2 = 1$ and this brings in a cost of $\alpha_2$. A similar approach can be used to consider the case with  $c_1 = 1$ (\Cref{star-base-cond-c1-1}) where the two possibilities are $t =1$ and $t=2$. If $ t= 1$, it must be that $c_2 = 0$ and therefore the cost incurred is only $\alpha_1$ while if $t = 2$, then it must be that $c_2 = 1$ and therefore the cost incurred is $\alpha_1 + \alpha_2 + \beta_{12}$.

Given these base conditions, we are now ready to compute the value of the function $g(i, \cdot, \cdot)$, for $i \geq 3$, in terms of the value corresponding to smaller sub-trees. For a given $t$,
\begin{align}
g(i, 0, t) &= \min( g(i-1, 0, t), g(i-1, 0, t - 1) + \alpha_i)\label{star-00}\\
g(i, 1, t) & = \min( g(i-1,1,t), g(i-1,1,t-1) + \alpha_i + \beta_{1,i}) \label{star-10}
\end{align}
The trees corresponding to $g(i, \cdot, \cdot)$ and  $g(i-1,  \cdot, \cdot)$ are shown in \Cref{fig:star-graph-5} for $i=4$.
When $t$ nodes are to be chosen from $\{1, \ldots, i \}$ with $c_1 = 0$, $c_i$ can take values either $0$ or $1$. If $c_i = 0$,  we must select $t$ nodes from $\{1, \ldots, i-1 \}$, while if $c_i =1$, we must select $t-1$ nodes from $\{1, \ldots, i-1 \}$. If $c_i=0$, there is no additional cost incurred while when $c_i = 1$, an additional cost of $\alpha_i$ is incurred. These two cases give rise to \Cref{star-00}. Note that the range of permissible value of $t$ varies for these two cases and can be similar identified as in the recursions in the case of the series graph. For example, for the case $c_1=0$ and $c_{i} =  0$, $t$ can only range from $0$ to $i-2$ in $g(i,0,t)$, while for $c_1=0$ and $c_{i} =1$, $t$ can range from  $1$ to $i-1$.
Using similar logic,  \Cref{star-10} can be written for the case where $c_1 = 1$. Finally, $Q_{star}$ is obtained by looking at all possible values of $g(n, \cdot, \cdot)$ as follows:
\begin{align}
Q_{star} &= \displaystyle \min\left(\min_{ t_1 \in [k, n-1]} g(n, 0, t_1), \min_{t_2 \in [k, n] }g(n, 1, t_2)\right) \label{star-root-0}
\end{align}
This gives rise to the linear programming formulation by linearization as in Proposition \ref{prop:series_graph_lp_ccqkp}.
\hfill \halmos
\subsection{Formulation for General Trees}
\label{sec:CCQKP-Tree}
We now provide the solution to the problem on general trees. A tree graph has several star graphs and series graphs as its components. The algorithm for a series graph involved a bottom-up traversal from the leaf node to the root while the star graph algorithm implicitly involved a traversal from the left most node to the right most node. The dynamic programing algorithm will involve solving sub-problems on a left to right as well as bottom up traversal of the nodes of the tree. Given a tree $T = (V,E)$, we are particularly interested in,
\begin{align}
\label{opt:card-constr-quad-knapsack-trees}
Q_{T} = \min \left\lbrace \sum_{i=1}^n \alpha_i c_i + \sum_{(i,j) \in E} \beta_{i,j} c_i c_{j}:  \sum_{i=1}^n c_i \geq k, c_i \in \{0,1\} \forall i \in [n] \right\rbrace
\end{align}

\begin{theorem} \label{thm:dp-solution}
The value of $Q_{T}$  in \eqref{opt:card-constr-quad-knapsack-trees} can be obtained using the following linear program:
\begin{align}
Q_{T} = & \max_{z, \bold{x}} z  \nonumber \\
\text{s.t } &   z \leq  x_{1, d_{1}, 0, t},\,\,\, t \in [k, n-1] \nonumber \\
& z \leq x_{1, d_{1}, 1, t},\,\,\, t \in [k, n] \nonumber \\
& x_{i, s, 0, 0} = 0\; \text{for }i \in V(i,s), \; s = [0,  d_i ] \nonumber  \\
& x_{i, s, 1, 1 } = \alpha_i \; \text{for } i  \in V(i,s), \;s = [0,  d_i] \nonumber  \\
&\text{For each internal node } i: \nonumber\\
& \;\;\; x_{i, 1, 0, t} \leq x_{i(1), d_{i(1)}, 0, t} \;\; \text{for } t=[0,  N(i,1)-2] \label{internal-0-0-0}\\
& \;\;\; x_{i, 1, 0, t} \leq x_{i(1), d_{i(1)}, 1, t} \;\; \text{for }t=[1, N(i,1)-1] \label{internal-0-1-0}\\
& \;\;\; x_{i, 1, 1, t} \leq x_{i(1), d_{i(1)}, 0, t-1} + \alpha_i \;\; \text{for }t=[1,  N(i,1)-1 ] \label{internal-1-0-0}\\
& \;\;\;  x_{i, 1, 1, t} \leq x_{i(1), d_{i(1)}, 1, t-1} + \alpha_i + \beta_{i, i(1)} \;\; \text{for } t=[2,  N(i,1)] \label{internal-1-1-0}\\
& \text {For each internal node i with out-degree at least 2:} \nonumber \\
& \;\;\; x_{i, s, 0, t} \leq x_{i, s-1, 0, t-a} + x_{i(s), d_{i(s)}, 0, a}  \;\;\text{for } t =[0, N(i,s)-2],  a \in [a_{\text{min}}^1 , a_{\text{max}}^1]  \label{internal-s-0-0}\\
 & \;\;\;  x_{i, s, 0, t} \leq x_{i, s-1, 0, t-a} + x_{i(s), d_{i(s)}, 1, a}, \;\;\text{for } t =[1,  N(i,s)-1], a \in [a_{\text{min}}^2 , a_{\text{max}}^2]  \label{internal-s-0-1}\\
  & \;\;\;  x_{i, s, 1, t} \leq x_{i, s-1, 1, t-a} + x_{i(s), d_{i(s)}, 0, a}, \;\;\text{for } t =[1, N(i,s)-1], a \in [a_{\text{min}}^3 , a_{\text{max}}^3 ] \label{internal-s-1-0}\\
   & \;\;\;  x_{i, s, 1, t} \leq x_{i, s-1, 1, t-a} + x_{i(s), d_{i(s)}, 1, a} + \beta_{i, i(s)}, \;\;\text{for } t =[2, N(i,s)], a \in [a_{\text{min}}^4 , a_{\text{max}}^4 ] \label{internal-s-1-1}
\end{align}
where,\\
 $a_{min}^1 = \max(0, t- (N(i,s-1) -1)), a_{max}^1 = \min( N(i(s), d_{i(s)})-1,t) $,\\
$a_{min}^2 = \max(1, t- (N(i,s-1) -1)), a_{max}^2 = \min( N(i(s), d_{i(s)}),t)$,\\
$a_{min}^3 = \max(0, t- (N(i,s-1))), a_{max}^3 = \min( N(i(s), d_{i(s)}-1),t-1)$,\\
$a_{min}^4 = \max(1, t- (N(i,s-1))), a_{max}^4 = \min( N(i(s), d_{i(s)}),t-1)$.
\end{theorem}
\textit{Proof:} Assume the tree is endowed with a fixed ordering of the nodes. For any node $i$, we denote by $x_{i, s, y, t}$  the minimum value of the sub-problem where: (1) the vertices and edges are restricted to the tree $T(i,s)$ rooted at $i$ and containing the first $s$ sub-trees of $i$ (as per the given ordering), (2) $c_i$ takes a value $y \in \{0,1\}$ and (3) exactly $t$ nodes in $T(i,s)$ take a value of $1$. This is given as:
\begin{align}
\label{opt:dp-defn-f}
x_{i, s, y, t}= \min \left\lbrace \sum\limits_{l \in V(i,s)} \alpha_l c_l + \sum\limits_{(l,j) \in V(i,s)} \beta_{lj} c_l c_j:  \sum_{l \in V(i,s)} c_l = t, c_i = y, c_l \in \{0,1\} \forall l \in V(i,s)  \right\rbrace
\end{align}
where $s \in \{0, \ldots, d_i \} $, $y \in \{0,1\}$, $t \in \{ y, \ldots, N(i,s)-1 + y \}$. Then the optimal value of the overall problem (\ref{opt:card-constr-quad-knapsack-trees}) can be expressed in terms of the sub-problems as follows:
\begin{align}
Q_{T} & \leq   x_{1, d_{1}, 0, t}, \,\, t \in [k, n-1] \\
 Q_{T} & \leq x_{1, d_{1}, 1, t} \,\, t \in [k, n]
\end{align}
Note the similarity of the above updates with the updates for $Q_{star}$ on star graphs.  The function $g(\cdot, \cdot, \cdot)$ previously considered, looked at the sub-tree rooted at node 1 always. The first argument to $g(\cdot ,\cdot, \cdot)$ considered $n-1$ values in $[2,n]$ corresponding to the children of node $1$. The second index in  $x_{1 ,\cdot, \cdot, \cdot}$ captures similar information. \\
We will now discuss the recursions to compute $\bold{x}$. The base cases include similar base cases as the star graph and the series case as follows:
\begin{align*}
x_{i, s, 0, 0} &= 0, \,\,\, i \in V(i,s), \,\, s \in [ 0,  d_i]  \\
x_{i, s, 1, 1 } &= \alpha_i , \,\,\, i  \in V(i,s), \,\,  s \in [ 0, d_i]
\end{align*}
The value $x_{i,s,0,0}$ looks at the case where $c_i=0$ and the entire sub-tree  rooted at $i$ containing the first $s$ subtrees of $i$ have $0$ nodes selected. This automatically suggests that for every $j$ in $V(i,s)$, $c_j=0$. In this case the optimal value of the sub-problem is 0.  On the other hand,  $x_{i, s, 1, 1 }$ deals with the case where $c_i = 1$ and exactly one of the nodes in $T(i,s)$ is set to $1$ (and this node has to be $c_i$ by definition). Therefore the optimal value is $\alpha_i$. This is true for every valid value of $s$.

Next, we will consider the recursions for $x_{i, 1, \cdot, \cdot}$ which only looks at the tree rooted at $i$ containing all nodes below and including the first child of $i$. The recursions are provided by \Cref{internal-0-0-0,internal-0-1-0,internal-1-0-0,internal-1-1-0} and follow from the underlying recursions  in the series graph. For  an internal node $i$, let $j=i(1)$ denote its first child.  Since $c_j \in \{0,1\}$, $x_{i, 1, 0, t}$ must take the minimum value out of $x_{j, d_j, 0, t}$ and $x_{j, d_j, 1, t}$ for all feasible values of $t$, while $x_{i,1,1,t}$ must take the minimum value from $x_{j, d_j, 0, t-1}$ and $x_{j, d_j, 1, t-1}$ . These recursions are same as the recursions in \Cref{rec:series1,rec:series3} for $\bar{f}(\cdot, \cdot, \cdot)$ with $i+1$ replaced with the first child $i(1)$.

\begin{figure}[h]
\centering
\begin{tikzpicture}[edge from parent/.style={draw,-latex}]
 \tikzstyle{level 1}=[sibling distance=20mm]
  \tikzstyle{level 2}=[sibling distance=12mm]
\node (1) [] {i}
    child { node (2) [label={above left: T1}]{ i(1)}
      child { node (5)  [label={below right:$ x_{i,2,\cdot,\cdot}$}]   {}
    }
      child { node (6) {}
    }
  }
    child { node (3) []{i(2)}
      child { node (7){}
    }
  }
  child { node (4)[label={above right: T2}]{i(3)}
      child { node (8){}
    }
     child { node (10)[]{}
    }
     child { node (9)[label={below right:$ x_{i(3), d_{i(3)},\cdot,\cdot}$}]{}
    }
  };

  \begin{pgfonlayer}{background}

  \draw[red,fill=gray,opacity=0.15](1.north east) to[closed,curve through={(3.east)..  (7.south east).. (6.south east).. (5.south west) .. (2.west) }] (1.north);

    \draw[blue,fill=gray,opacity=0.6](4.north east) to[closed,curve through={(9.north east).. (9.south west) ..  (8.west) }] (4.north);
  \end{pgfonlayer}
\end{tikzpicture}
\caption{$x_{i,3,\cdot,t}$ can be computed using $x_{i, 2, \cdot, t-a}$ (corresponding to the sub-tree T1) and   $x_{i(3), d_{i(3)}, \cdot, a}$ (corresponding to the sub-tree T2) so that $t-a$ nodes are selected from T1 and $a$ nodes are selected from T2.}
\label{fig:tree-computation}
\end{figure}
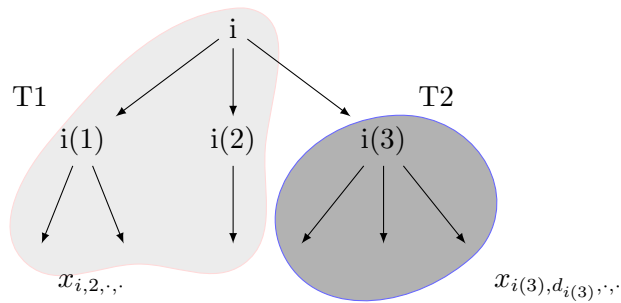

Now for an internal node $i$ with at least $2$ children, we compute $x_{i, s, \cdot, \cdot}$ using $x_{j, \cdot, \cdot, \cdot}$ corresponding to all the children $j$ of $i$, as well as $x_{i, s-1, \cdot, \cdot}$. The two sub-trees involved in computing $x_{i, s, \cdot, \cdot}$  are depicted as shaded regions (labelled T1 and T2) in \Cref{fig:tree-computation}, for $s=3$. First, suppose $c_i = 0$. We are interested to compute $x_{i, s, 0, t}$. The $t$ nodes to be selected from the tree rooted at node $i$ can be split  between $T_1$ and $T_2$ in various ways. The value of $t$ itself can only range between $0$ to $N(i, s) - 1$ since $c_i = 0$. Suppose $a$ nodes are selected in $T_2$ and $t-a$ nodes are selected in $T_1$. The range of $a$ differs based on whether the root node of $T_2$ is selected or not.  Suppose the root node of $T_2$  is not selected (that is, $c_{i(s)} = 0$). Since the number of nodes selected from $T_2$ is $a$ and $i(s)$ is not selected,
\begin{align}
 0 \leq a \leq N_{T_2}-1
 \end{align} where $N_{T_2}  $ is the number of nodes in $T_2$. Also  since the number of nodes selected from $T_1$ is $t-a$ and $i$ is not selected,
 \begin{align}
  0 \leq t-a  \leq N_{T_1} - 1
  \end{align} Based on the above two inequalities, we get $a \in [a_{min}^1 , a_{max}^1]$,
  where $a_{min}^1 = \max(0, t - (N_{T_1} - 1))$ and $a_{max}^1 = \min(N_{T_2} - 1, t )$.
  For these values of $a$, the cost is just $x_{i, s-1, 0, t-a} +x_{i(s), d_{i(s)},0, a}$ and no additional cost gets added as both $i$ and $i(s)$ are not selected. Hence we get \Cref{internal-s-0-0}. A similar treatment gives us \Cref{internal-s-0-1} corresponding to $c_{i(s)}=1$. Note that the cost $\alpha_{i(s)}$ is already part of $x_{i(s), d_{i(s)}, 0, a}$ and does not need to be explicitly added. $x_{i, s, 0, t}$ must take the minimum value of all terms in the RHS of  \Cref{internal-s-0-0,internal-s-0-1}.

Now suppose $c_i = 1$, we are interested to compute $x_{i, s , 1, t}$.  The valid values of $t$ range from $1$ to $N(i, s)$. Here we will illustrate the case where $c_{i(s)} = 1$.  Again assume the sub-tree $T_2$ rooted at $i(s)$ contains $a$ nodes and $T_1$ contains $t-a$ nodes. Using similar reasoning from the earlier step, the range of valid values for $a$  can be derived. Now in addition to the sub-tree costs $x_{i, s-1, 1, t-a} + x_{i(s), d_{i(s)}, 1, a}$, we also incur an additional cost $\beta_{i, i(s)}$ of selecting both nodes $i$ as well as $i(s)$. Hence we get \Cref{internal-s-1-1}. The case where $c_{i(s)}=0$ follows using similar logic(see \eqref{internal-s-1-0}). Note that the individual item costs $\alpha_i$ and $\alpha_{i(s)}$ are absorbed in $x_{i, s-1, 1, t-a}$ and $x_{i(s), d_{i(s)}, 1, a}$ respectively.  Finally, $x_{i, s, 1, t}$ must take the minimum value of all terms in the RHS of  \Cref{internal-s-1-0,internal-s-1-1} and hence the inequalities arise.
The variable $z$ denotes $Q_{T}$.
\hfill \halmos

 The particular sub-trees used in these cases illustrated in \Cref{fig:tree-computation} can be viewed as generalizations of the  sub-trees used in the star graph.  For example, in \Cref{fig:star-graph-5}, the darker sub-tree (with nodes $\{1, 2, 3 \}$) represents $T1$ and the sub-tree $T_2$ (not shown explicitly) trivially has exactly $1$ node (node 4 in the particular instance in the \Cref{fig:star-graph-5}).

The number of variables in the optimization is $\sum_{i=1}^n 2 d_i n $ or $O(n^2)$ (as $\sum_{i} d_i = O(n)$ for a tree) while the number of constraints is $O(n^3)$.

\section{Probability Bounds with Trees}
\label{sec:prob_bounds_tree}
We will now use the results from the previous section to derive tight bounds for the probability $U(k)$ when the univariate and bivariate probabilities for a given tree graph $T = (V, E)$ are known. The exponential sized dual linear program involves two sets of constraints \eqref{constr:1} and \eqref{constr:0} to compute $U(k)$. The results in \Cref{sec:ccqkp} discuss the separation problem involving  \eqref{constr:1}. We will now see that  \eqref{constr:0} can be equivalently represented by a set of linear constraints.

\begin{lemma}
\label{lemma:separation-problem-bqp-reformulation}
For a given set of values for $\{ \lambda, \boldsymbol{\alpha}, \boldsymbol{\beta} \}$, constraint (\ref{constr:0}) is equivalent to verifying the feasibility of the following set of linear constraints in the variables  $\left\lbrace  \boldsymbol{\eta}, \boldsymbol{\gamma}, \boldsymbol{\Delta}, \boldsymbol{\tau}, \boldsymbol{\chi}\right\rbrace$ :
\begin{align}
&\lambda - \sum\limits_{(i,j) \in E} \left( \Delta_{ij} + \chi_{ij} \right) - \sum\limits_{i=1}^n \tau_i \geq 0 \label{sep-prob:const0}\\
 & \sum\limits_{j: (i,j) \in E} \left( \Delta_{ij} - \eta_{ij} \right) + \sum\limits_{j: (j,i) \in E} \left( \Delta_{ji} - \gamma_{ji} \right) + \tau_i + \alpha_i  \geq 0 \;\; \text{ for }i \in [n]\label{sep-prob:const1}\\
& \eta_{ij} + \gamma_{ij} - \Delta_{ij} + \chi_{ij} + \beta_{ij} \geq 0 \;\; \text{ for }(i,j) \in E  \label{sep-prob:const2}\\
& \boldsymbol{\eta}, \boldsymbol{\gamma}, \boldsymbol{\Delta}, \boldsymbol{\tau}, \boldsymbol{\chi} \geq 0 \label{sep-prob:const3}
\end{align}
\end{lemma}

\textit{Proof:} Constraint \eqref{constr:0} may be re-written as,
\begin{align}
\label{separation:BQP}
\lambda  - \max\limits_{\bold{c} \in  \{0,1\}^n} \left\lbrace   -\sum_{i=1}^n \alpha_i c_i - \sum_{(i,j) \in E} \beta_{ij} c_i c_j:   c_i \in \{0,1\} \forall \; i \right\rbrace \geq 0
\end{align}
Note that for a given set of values for $\{\lambda, \boldsymbol{\alpha},  \boldsymbol{\beta}\}$, this reduces to optimizing a quadratic function over the extreme points of the unit hypercube which corresponds to optimization over the Boolean quadric polytope.
We can now directly apply the results in \cite{Padberg1989} who derived a tight formulation when the graph is a tree. The results therein implies the following LP relaxation is tight for the maximization problem in Eqn (\ref{separation:BQP}) when the sparsity pattern of the quadratic terms $\boldsymbol{\beta}$ is given as tree  $T=(V,E)$:
\begin{align*}
\max\limits_{\bold{c}, \bold{y}} &  -\sum_{i=1}^n \alpha_i c_i - \sum_{(i,j) \in T} \beta_{ij} y_{ij}\\
\textrm{s.t } & y_{ij} - c_i \leq 0 \;  \text{ for } (i,j) \in E \\
&  y_{ij} - c_j \leq 0 \; \text{ for } (i,j) \in E \\
& c_i + c_j - y_{ij} \leq 1 \; \text{ for } (i,j) \in E\\
& 0 \leq c_i \leq 1 \; \text{ for } i \in [n] \\
& 0 \leq y_{ij} \leq 1 \; \text{ for } (i,j) \in E
\end{align*}
 The constraints  \eqref{sep-prob:const1}, \eqref{sep-prob:const2} and \eqref{sep-prob:const3} in Lemma \ref{lemma:separation-problem-bqp-reformulation} correspond to the constraints in the dual of the above linear optimization problem. Constraint \eqref{sep-prob:const0} appears as a result of  forcing the dual objective to be non-negative (due to the non-negativity requirement of the objective in \eqref{separation:BQP}).
 \hfill \halmos

 We are now ready to provide the main result of the paper.

\begin{theorem} \label{bivar-tree-theorem} Consider an $n$-dimensional Bernoulli random vector $\tilde{\bold{c}}$ where the univariate probabilities $p_i = P(\tilde{c}_i = 1) \text{ for } i = 1, \ldots, n$ and the bivariate marginals $p_{ij} = P(\tilde{c}_i = 1, \tilde{c}_j = 1)$ for all $(i,j) \in E$ for a tree graph are specified.
Let $U^*_k$ denote the optimal value of the following linear program over the decision variables $\Gamma = \{ \lambda, \boldsymbol{\alpha}, \boldsymbol{\beta}, \boldsymbol{\Delta}, \boldsymbol{\eta}, \boldsymbol{\chi}, \boldsymbol{\tau}, \boldsymbol{x}\}$:
\begin{align*}
U^*_k= \displaystyle \min_{\Gamma}\; & \lambda + \sum\limits_{i=1}^n \alpha_i p_i + \sum\limits_{(i,j) \in E} \beta_{ij} p_{ij}\\
\textrm{s.t.} & \lambda - \sum\limits_{(i,j) \in E} \left(\Delta_{ij} + \chi_{ij} \right) - \sum\limits_{i=1}^n \tau_i \geq 0 \\
 & \sum\limits_{j: (i,j) \in E} \left( \Delta_{ij} - \eta_{ij} \right) + \sum\limits_{j: (j,i) \in E} \left( \Delta_{ji} - \gamma_{ji} \right) + \tau_i + \alpha_i  \geq 0 \;\; \text{for } i \in [n] \\
& \eta_{ij} + \gamma_{ij} - \Delta_{ij} + \chi_{ij} + \beta_{ij} \geq 0 \;\; \text{for } (i,j) \in E \\
& \lambda + z \geq 1  \\
& x_{1, d_1, 0, t} - z \geq 0 \;\; \text{for } t \in [k,n]  \\
& x_{1, d_1, 1, t} - z \geq 0 \;\; \text{for } t \in [k,n]  \\
& x_{i,s,0,0} = 0 \;\; \text{ for }i \in [n], \; \text{for } s \in [0, d_i] \\
& x_{i,s,1,1} - \alpha_i = 0 \;\; \text{for } i \in [n], \;  \text{for } s \in [0, d_i] \\
& \text {For each internal node i: }  \\
& \;\;\;\;\;\;  x_{i(1), d_{i(1)}, 0, t} - x_{i,1,0,t} \geq 0 \;\; \text{for } t \in [0, N(i,1)-2] \\
&  \;\;\;\;\;\; x_{i(1), d_{i(1)}, 1, t} - x_{i,1,0,t} \geq 0 \; \; \text{for } t \in [1,  N(i,1)-1]  \\
&  \;\;\;\;\;\; x_{i(1), d_{i(1)}, 0, t-1} - x_{i,1,1,t} + \alpha_i \geq 0 \;\;\text{for }t \in [1, N(i,1)-1] \\
&  \;\;\;\;\;\;  x_{i(1), d_{i(1)}, 1, t-1} - x_{i,1,1,t} + \alpha_i  + \beta_{i, i(i)} \geq 0 \;\; \text{for } t \in [2, N(i,1)] \\
& \text {For each internal node i with out-degree at least 2:} \\
& \;\;\;\;\;\; x_{i, s-1, 0, t-a} + x_{i(s), d_{i(s)}, 0, a} - x_{i,s,0,t} \geq 0
  \text{ for }s =[ 2,  d_i], t= [0 , (N(i,s)-2)], a= [a_{min}^1, a_{max}^1] \\
& \;\;\;\;\;\; x_{i, s-1, 0, t-a} + x_{i(s), d_{i(s)}, 1, a} - x_{i,s,0,t} \geq 0
 \text{ for } s = [2, d_i], t= [1 , (N(i,s)-1)], a= [a_{min}^2 , a_{max}^2] \\
& \;\;\;\;\;\; x_{i, s-1, 1, t-a} + x_{i(s), d_{i(s)}, 0, a} - x_{i,s,1,t} \geq 0
 \text{ for } s = [2,  d_i], t= [1, (N(i,s)-1)], a= [a_{min}^3, a_{max}^3]\\
& \;\;\;\;\;\; x_{i, s-1, 1, t-a} + x_{i(s), d_{i(s)}, 1, a} - x_{i,s,1,t} + \beta_{i,i(s)}\geq 0
\text{ for } s = [2,  d_i], t= [2 ,(N(i,s))], a= [a_{min}^4, a_{max}^4] \\
  & \boldsymbol{\eta}, \boldsymbol{\gamma}, \boldsymbol{\Delta}, \boldsymbol{\tau}, \boldsymbol{\chi} \geq 0
\end{align*}
where $a_{min}^1 = \max(0, t- (N(i,s-1) -1)), a_{max}^1 = \min( N(i(s), d_{i(s)}-1),t) $,\\
$a_{min}^2 = \max(1, t- (N(i,s-1) -1)), a_{max}^2 = \min( N(i(s), d_{i(s)},t)$,\\
$a_{min}^3 = \max(0, t- (N(i,s-1))), a_{max}^3 = \min( N(i(s), d_{i(s)}-1),t-1)$,\\
$a_{min}^4 = \max(1, t- (N(i,s-1))), a_{max}^4 = \min( N(i(s), d_{i(s)}),t-1).$\\
Then, $U(k) = U^*_k$.
\end{theorem}
\textit{Proof:}
 We derive the linear programming reformulation  by considering each of the two groups of constraints  \eqref{constr:1} and \eqref{constr:0}  that arise in the dual problem. Constraint \eqref{constr:1} can be re-written  as,
\begin{align*}
 \lambda  + \min\limits_{\bold{c} \in \{0,1\}^n} \left\lbrace \sum_{i=1}^n \alpha_i c_i + \sum_{(i,j) \in T} \beta_{ij} c_i c_j:  \sum_{i=1}^n c_i \geq k \right\rbrace \geq 1
\end{align*}
This is the cardinality constrained quadratic knapsack problem in minimization form.  In \Cref{thm:dp-solution}, we provided a linear programming reformulation for this sub-problem. Plugging in the constraints from the linear program in \Cref{thm:dp-solution} and forcing the objective value to be greater than $1$ gives us all, except the first three  constraints in the linear program. The first three constraints are equivalent to \Cref{constr:0} as stated in \Cref{lemma:separation-problem-bqp-reformulation}.
This completes the proof of \Cref{bivar-tree-theorem}. \hfill \halmos

The number of variables and constraints in our linear programming solution are  $O(n^2)$ and $O(n^3)$ respectively, which implies that the tight bound is solvable in polynomial time.

\subsection{Probability with a Tree Graphical Model}
\label{sec:cond-ind}
We now relate our bounds to the computation of $\mathbb{P}(\sum_{i=1}^n \tilde{c}_i \geq k)$ under the same information as before but now focus on a tree graphical model where conditional independence is assumed. The difference from the bounds in the paper is that this induces a unique distribution where every random variable $\tilde{c}_i$ is independent of all its  siblings, conditional on knowledge of the realization of its parent $i$. The next proposition provides a dynamic programming recursion (similar to \Cref{sec:ccqkp}) to compute $\mathbb{P}(\sum_{i=1}^n \tilde{c}_i \geq k)$ in this case.

\begin{proposition}
\label{thm:cond_ind_recursion} Consider a $n$ dimensional Bernoulli random vector $\tilde{\bold{c}}$ where the univariate probabilities $p_i = P(\tilde{c}_i = 1) \text{ for } i = 1, \ldots, n$ and the bivariate marginals $p_{ij} = P(\tilde{c}_i = 1, \tilde{c}_j = 1)$ for all $(i,j) \in E$ for a tree graph are specified and assume it is a conditionally independent distribution on the tree as in (\ref{eqn:cond-ind-whole-tree}). Then,
the following recursions can be used to compute $\mathbb{P}(\sum_{i=1}^n \tilde{c}_i \geq k)$:
\begin{align}
&w_{i,s,0,0} = p_0(i) \,\,\text{ for } i \in [n], s \in [0, d_i] \\
& w_{i,s,1,1}  =p_1(i) \,\,\text{ for } i \in [n], s \in [0, d_i] \\
& \text{For each internal node i}, \nonumber\\
 & w_{i, 1, 0, t} = \frac{w_{i(1), d_{i(1)}, 0, t} \times p^{00}(i,i(1)) \times\mathbbm{1}\{t\in [0, N(i,1)-1]\}}{p_0(i(1))}  \nonumber\\
 &\,\,\,\,\,\,\,\, + \frac{w_{i(1), d_{i(1)}, 1, t} \times p^{01}(i,i(1)) \times \mathbbm{1}\{ t \in[1, N(i,1)-1]\}}{p_1(i(1))} ,\,\,\,  t \in [0, N(i,1)-1] \label{rec:w-1-0-t}\\
  & w_{i, 1, 1, t } = \frac{w_{i(1), d_{i(1)}, 0, t-1} \times p^{10}(i,i(1)) \times \mathbbm{1}\{t \in [1, N(i,1)-1]\}}{p_0(i(1))}  \nonumber \\
 &\,\,\,\,\,\,\,\,  + \frac{w_{i(1), d_{i(1)}, 1, t-1}\times p^{11}(i,i(1)) \times \mathbbm{1}\{ t \in[2, N(i,1)]\}}{p_1(i(1))} ,\,\,\,  t \in [1, N(i,1)-1]  \label{rec:w-1-1-t}\\
& \text {For each internal node i with out-degree at least 2:} \nonumber \\
& \;\;\;\;\;\;  w_{i,s,0,t} = \sum_{a = a_{min}^1}^{a_{max}^1}\frac{ w_{i, s-1, 0, t-a} \times w_{i(s), d_{i(s)}, 0, a} \times p^{00}(i,i(s))\mathbbm{1}\{ t \in [0, (N(i,s)-2)] \}}{p_0(i) p_0(i(s))} \nonumber \\
& \;\;\;\;\;\;\;\;\;\;  + \sum_{a = a_{min}^2}^{a_{max}^2} \frac{ w_{i, s-1, 0, t-a} \times w_{i(s), d_{i(s)}, 1, a} \times p^{01}(i, i(s)) \mathbbm{1}\{ t \in [1, (N(i,s)-1)] \}} {p_0(i) p_1(i(s))} \nonumber \\
&  \;\;\;\;\;\;\;\;\;\;  \;\;\;\;\;\;\;\;\;\;  \;\;\;\;\;\;\;\;\;\; \;\;\;\;\;\;\;\;\;\; \,\,\text{ for } s \in [2,  d_i] , t \in [ 0, N(i,s) - 1 ] \label{rec:w-s-0-t} \\
& \;\;\;\;\;\;  w_{i,s,1,t} = \sum_{a = a_{min}^3}^{a_{max}^3}\frac{ w_{i, s-1, 1, t-a} \times w_{i(s), d_{i(s)}, 0, a}
\times p^{10}(i,i(s))\mathbbm{1}\{ t \in [1, (N(i,s)-1)] \}}{p_1(i) p_0(i(s))} \nonumber \\
& \;\;\;\;\;\;\;\;\;\;  + \sum_{a = a_{min}^4}^{a_{max}^4} \frac{ w_{i, s-1, 1, t-a} \times w_{i(s), d_{i(s)}, 1, a} \times p^{11}(i, i(s)) \mathbbm{1}\{ t \in [2, N(i,s)] \}} {p_1(i) p_1(i(s))}   \nonumber \\
&  \;\;\;\;\;\;\;\;\;\;  \;\;\;\;\;\;\;\;\;\;  \;\;\;\;\;\;\;\;\;\; \;\;\;\;\;\;\;\;\;\; \,\,\text{ for } s \in [2,  d_i] , t \in [1,  N(i,s)] \label{rec:w-s-1-t}
\end{align}
where $a_{min}^1 = \max(0, t- (N(i,s-1) -1)), a_{max}^1 = \min( N(i(s), d_{i(s)}-1),t) $,\\
$a_{min}^2 = \max(1, t- (N(i,s-1) -1)), a_{max}^2 = \min( N(i(s), d_{i(s)}),t)$,\\
$a_{min}^3 = \max(0, t- (N(i,s-1))), a_{max}^3 = \min( N(i(s), d_{i(s)}-1),t-1)$,\\
$a_{min}^4 = \max(1, t- (N(i,s-1))), a_{max}^4 = \min( N(i(s), d_{i(s)}),t-1)$,\\
$p^{00}(i,j) = \mathbb{P}(\tilde{c}_i=0, \tilde{c}_j=0) = 1- p_i - p_j + p_{ij} ,\\
p^{01}(i,j) = \mathbb{P}(\tilde{c}_i=0, \tilde{c}_j=1) = p_j - p_{ij}, \\
p^{10}(i,j) = \mathbb{P}(\tilde{c}_i=1, \tilde{c}_j=0) = p_i - p_{ij} $, $p^{11}(i,j) =p_{ij} $
and $p_1(i) = p_i, p_0(i) = 1 - p_i$.

In particular, $\mathbb{P}(\sum_{i=1}^n \tilde{c}_i \geq k) = \sum_{t=k}^n w_{1,d_{1},0,t} +w_{1,d_{1},1,t}$.
\end{proposition}
The  proof idea is to express the probability $ \mathbb{P}\left(\sum_{l \in T(i,s)}\tilde{c}_l=  t, \tilde{c}_i = y \right)$, denoted by $w_{i, s, y, t}$, in terms of probabilities for smaller sub-trees of $T(i,s)$.  The details are provided in the online companion.

The approach for deriving the recursions  on $\bold{w}$ progresses in a very similar manner to the recursions for $\bold{x}$ in the quadratic knapsack problem on tree graphs in \Cref{sec:CCQKP-Tree}. Both approaches work on $O(n^2)$ variables.  The same sub-trees (as depicted in \Cref{fig:tree-computation})  are used in the computation of $x_{i,s,y,t}$ and $w_{i,s,y,t}$ in Theorems \ref{thm:dp-solution} and \ref{thm:cond_ind_recursion} respectively. For example, in \eqref{rec:w-s-1-t},the computation of  $w_{i,s,1,t}$ makes use of sums and products involving $w_{i,s-1,0,t-a}, w_{i(s), d_{i(s)},0,a}$ and $w_{i(s), d_{i(s)},1,a}$. In a similar vein, the computation of  $x_{i,s,1,t}$ in \Cref{bivar-tree-theorem} makes use of inequalities and sums  involving $x_{i,s-1,0,t-a}, x_{i(s), d_{i(s)},0,a}$ and $x_{i(s), d_{i(s)},1,a}$. For the sub-problems involving only the first sub-tree in \eqref{rec:w-1-0-t},  $w_{i,1,0,t}$ makes use of summations involving $w_{i(1), d_{i(1)}, 0,t}$ and $w_{i(1), d_{i(1)}, 1,t}$. Similarly the computation of $x_{i,1,0,t}$ makes use of inequalities involving  $x_{i(1), d_{i(1)}, 0,t}$ and $x_{i(1), d_{i(1)}, 1,t}$. Of course, the expressions themselves are different in the two theorems - Theorem \ref{thm:dp-solution} looks at solving an underlying optimization problem while the goal of  \Cref{thm:cond_ind_recursion}  is to simplify the computation of a probability which otherwise involves an exponential number of operations. A comparison of the underlying techniques proposed for these settings is provided in \Cref{tab:comparison-ub-ci}.

\begin{table}[h!]
\centering
\caption{Comparison of proposed approach for computing the upper bound $\max \mathbb{P}(\sum_{i=1}^k \tilde{c}_i \geq k)$ vs $ \mathbb{P}(\sum_{i=1}^k \tilde{c}_i \geq k)$ for the conditionally independent distribution on a given tree }
\label{tab:comparison-ub-ci}
\begin{tabular}{|p{3.5cm}|p{6cm}|p{6cm}|}
\hline
& \textbf{Upper Bound} & \textbf{Conditional Independence} \\\hline
Formulation & Linear program & Dynamic programming recursion\\\hline
No. of variables & $O(n^2)$ & $O(n^2)$ \\\hline
No. of operations/constraints& $O(n^3)$ constraints in the formulation &  $O(n^3)$ summations in \Cref{rec:w-1-0-t,rec:w-1-1-t,rec:w-s-0-t,rec:w-s-1-t}\\\hline
Type of operations & Inequalities involving summations over the variables & Equalities involving summations and products over the variables\\\hline

\end{tabular}

\end{table}
\section{Generalizations}
In this section, we discuss some generalizations where the proposed probability bounds can be applied.
\subsection{Lower Bound}
\label{sec:lower_bound}
 The analogous  approach can be used to find lower bounds $L(k)$ for the same sum on the given tree graph $T = (V,E)$. To see this, let $\tilde{d}_i = 1 - \tilde{c}_i$ and define $q_i = 1- p_i$ and $q_{ij} = \mathbb{P}(\tilde{c}_i=0, \tilde{c}_j = 0)$ for all $(i,j) \in E$. Then,
\begin{align}
\label{eq:lower_bound}
L(k) =  \min \mathbb{P}_{\theta}(\sum_{i=1}^n \tilde{c}_i  \geq k ) &= 1 + \min  - \mathbb{P}_{\theta}(\sum_{i=1}^n \tilde{c}_i  < k ) \\
 &= 1 - \max \mathbb{P}_{\theta}( \sum_{i=1}^n1 - \tilde{c}_i  \geq n - k + 1 ) \\
 & =  1 - \max \mathbb{P}_{\theta}( \sum_{i=1}^n \tilde{d}_i \geq n - k + 1) = U(n-k +1)
\end{align}
where $U(n-k +1)$ must be computed for the given tree graph  $T$ by setting $p_i = q_i$ and $p_{ij} = q_{ij}$ in \Cref{bivar-tree-theorem}.

\subsection{Bounds for Weighted Sum of Probabilities} Our approach can be generalized to compute upper and lower bounds for the weighted sum of probabilities $ \sum_{s=0}^{n} w_s \mathbb{P}(\sum_{i=1}^{n} \tilde{c}_i = s)$ with a given weight vector $\bold{w}  \in \mathbb{R}_{n+1}^{+}$. Weighted sums arise in a scenario where, for example, the set $\{1, \ldots, n\}$ can be partitioned into two  disjoint sets $A$ and $B$ such that the random variables $\{ \tilde{c}_i, i \in A \} $ are known to be mutually independent and also independent from $\{ \tilde{c}_j, j \in B \} $. The random variables corresponding to the set $B$ could however be dependent but the joint distribution over $\{ \tilde{c}_j, j \in B \} $ is unknown. Denote by $\tilde{\bold{c}}_A$ the vector of random variables corresponding to the set $A$. In particular, assume the following information:
\begin{enumerate}
\item A tree structure $T_B = (B, E_B)$
\item The univariate probabilities  $\mathbb{P}(\tilde{c}_i = 1) = p_i $, for $i \in [n]$
\item The bivariate probabilities $\mathbb{P}(\tilde{c}_i = 1, \tilde{c}_j = 1) = p_{ij} \text{ for } (i,j) \in E_B $
\item $\mathbb{P}(\tilde{c}_i = 1, \tilde{c}_j =1) = p_i p_j \text{ for }  i \in A, j \in B$
\item $\mathbb{P}(\tilde{\bold{c}}_A = \bold{r}_A) = \prod_{i \in A} \mathbb{P}(\tilde{c}_i = r_i )\text{ for }  \bold{r}_A \in \{0,1\}^{|A|}$
\end{enumerate}
By enumerating all ways in which $\sum_{i \in A} \tilde{c}_i$ and $\sum_{i \in B} \tilde{c}_i$ add up to a value $j \geq k$ we can express the tail probability as a weighted sum as below,
 \begin{align}
  \mathbb{P}(\sum_{i=1}^n \tilde{c}_i \geq k) &= \sum_{s=0}^{|B|}w_s \mathbb{P}(\sum_{i \in B} \tilde{c}_i = s)
   \label{eq:weighted_sum}
 \end{align}
 where  $w_s = \mathbb{P}(\sum_{i \in A} \tilde{c}_i \geq k - s) \geq 0 $. This relation follows as a consequence of independence between $\tilde{\bold{c}}_A$ and $\tilde{\bold{c}}_B$. $\sum_{i \in A} \tilde{c}_i$ is  a sum of independent but non-identical Bernoulli random variables and therefore takes a Poisson-binomial distribution, for which the  probability $w_s$ can be computed in a recursive manner in polynomial time (see \cite{chen1998}).

We now show that the computation of  $\max_{\theta \in \Theta} \sum_{s=0}^n w_s \mathbb{P}(\sum_{i=1}^n \tilde{c}_i = s)$, given any $\bold{w} \in \mathbb{R}_{n+1}^{+}$ can be easily done. The result in \Cref{bivar-tree-theorem} provides a linear program for the special  case where $w_s = 0$ for $s \in [0, k-1]$ and $w_s= 1$ for $s \in [k, n]$. Given any $\bold{w} \in \mathbb{R}_{n+1}^+$, the analogous exponential sized dual formulation is:
\begin{align}
\min\limits_{\lambda, \boldsymbol{\alpha},  \boldsymbol{\beta}}  &\; \lambda +  \sum_{i=1}^n \alpha_i p_i + \sum\limits_{(i,j) \in E} \beta_{ij} p_{ij} \nonumber \\
\textrm{s.t.} &
 \lambda + \sum_{i=1}^n \alpha_i c_i + \sum\limits_{(i,j) \in E} \beta_{ij} c_i c_j \geq w_s,   \;\; \text{ for } \bold{c} \in \{0,1\}^n \text{ where} \sum\limits_{i=1}^n c_i = s, \text{for } s \in [0,n] \label{constr:1_w} \\
& \lambda + \sum_{i=1}^n \alpha_i c_i + \sum\limits_{(i,j) \in E} \beta_{ij} c_i c_j \geq 0 \;\; \text{ for }  \bold{c} \in \{0,1\}^n  \label{constr:0_w}
\end{align}

Constraint (\ref{constr:0_w}) is exactly the same as constraint (\ref{constr:0}) and therefore \Cref{lemma:separation-problem-bqp-reformulation} gives a reformulation for this constraint. For constraint (\ref{constr:1_w}), for each value of $s$, a polynomial sized  linear programming formulation can be derived based on similar dynamic programming recursions in \Cref{sec:ccqkp}. The objective of the resulting linear program must be forced to take a value greater than $w_s$ (instead of a value of 1 in \Cref{bivar-tree-theorem}).

\subsection{Bounds for Order Statistics}
\label{sec:order-stats-bd}
We will now provide an application of our approach to obtain bounds on order statistics probabilities for random variables with any underlying distribution (that is, either discrete or continuous random variables).  Let $\tilde{X}= ( \tilde{X}_1, \ldots, \tilde{X}_n)$ denote $n$ real valued random variables. The order statistics is a re-arrangement of the $\tilde{X}_i$ denoted as $\tilde{X}_{1:n} \leq  \tilde{X}_{2:n}  \leq \ldots \tilde{X}_{k:n} \leq \ldots \leq \tilde{X}_{n:n}$. The $k^{th}$ order statistic $\tilde{X}_{k:n}$ denotes the $k^{th}$ smallest random variable among the $n$ random variables. The cases $k=1$ and $k=n$ corresponds to the minimum and maximum of the random variables respectively.

Computing the density function of the $k^{th}$ order statistic is a problem that has drawn the attention of researchers since several decades. The early  methods focussed on computing the probabilities of order statistics given the marginal distributional information of i.i.d random variables with subsequent extensions to the independent and non-identical random variables (see \citep{OrderStatsTextBook}).  A general technique for deriving the $k^{th}$ order statistic probabilities for dependent random variables has been provided in \cite{david2004order} but as pointed by the authors, simple expressions are usually possible only for specific cases or under more restrictions (e.g multivariate normal with equal and positive correlations in \cite{Tong1990}, exchangeable random variables in \cite{Arellano2007}, \cite{Arellano2008} etc).  There has also been interest in bounding the expected value of order statistics assuming moments of the random variables (see \cite{Rychlik1994,karthik2006,OrderStatsTextBook,david2004order}).

We will now show how our formulation can be used to compute bounds on the cdf of order statistics for dependent random variables. Specifically, we are interested in computing the probability $\mathbb{P}(\tilde{X}_{k:n} \leq x)$ for any given value $x$.  If the $k^{th}$ smallest random variable must take a value less than $x$, then it implies that at least $k$ of the random variables $\tilde{X}_1, \ldots, \tilde{X}_n$ must take a value less than $x$. Let $\tilde{c}_i = \mathbbm{1}\{ \tilde{X}_i \leq x \}$. Then,
\begin{align}
\label{eqn:orderstats-connection}
\mathbb{P}(\tilde{X}_{k:n} \leq x) = \mathbb{P}( \sum_{i=1}^n \tilde{c}_i \geq k)
\end{align}
Given a tree graph $T = (V, E)$ on $n$ nodes and a value $x \in \mathbb{R}$, let $\Theta_{x}$ denote the set of all joint distributions on $\bold{\tilde{X}}$ consistent with the univariate cdfs and bivariate cdfs evaluated at $x$ for random variables corresponding to edges in $T$,  as follows:
\begin{align}
\Theta_{x}= \big\{
  \theta: & \mathbb{P}_\theta\left(\tilde{X}_i \leq x, \tilde{X}_j \leq x
    \right) = p_{ij}(x) \text{ for }
    (i,j) \in E , \;\;  \mathbb{P}_\theta\left(\tilde{X}_i \leq x
    \right) = p_i(x) \text{ for }
    i \in [n]  \big\}.
    \end{align}

    We are interested in computing the largest possible and least possible values of the order statistics probabilities over all distributions in $\Theta_{x}$.
For any value $x \in \mathbb{R}$, define    $ U_{\text{os},k}(x)$ and $L_{\text{os},k}(x)$ as,
    \begin{align}
    U_{\text{os},k}(x)  &= \max_{\theta \in \Theta_{x}}  \mathbb{P}_{\theta}(\tilde{X}_{k:n} \leq x) \\
        L_{\text{os},k}(x)  &= \min_{\theta \in \Theta_{x}}  \mathbb{P}_{\theta}(\tilde{X}_{k:n} \leq x)
    \end{align}
The following proposition is an application of \Cref{bivar-tree-theorem}  and provides an upper bound for $  U_{\text{os},k}(x) $.
\begin{proposition}
\label{thm:order-stats-bd}
For any $x \in \mathbb{R}$, suppose we are given the following information:
\begin{enumerate}
\item A tree graph $T = (V,E)$ where $V = \{1, \ldots, n\}$ with node $1$ as the root node.
\item Univariate Probabilities $p_i(x) = P(\tilde{X}_i \leq x) \text{ for } i \in [n]$,
\item Bivariate probabilities $p_{ij}(x) = P(\tilde{X}_i \leq x, \tilde{X}_j \leq x)$ for all $(i,j) \in E$, consistent with the univariate probabilities $p_i(x)$.
\end{enumerate}

Let $\hat{U}_{\text{os},k}(x)$ denote the value of the linear program in \Cref{bivar-tree-theorem} supplied with the above information. Then $\hat{U}_{\text{os},k}(x) =  U_{\text{os},k}(x) $.
\end{proposition}
The proposition follows as a consequence of \eqref{eqn:orderstats-connection}. Note that $L_{\text{os},k}(x)$ can be computed as a function of  $U_{\text{os}, k}(n-k+1) $ with appropriate substitutions for the bivariate and univariate probabilities as described in \Cref{sec:lower_bound}.

\section{Numerical Computations}
We now present the results of our numerical computations. For convenience, in this section, we will use the notation $S = \sum_{i=1}^n \tilde{c}_i$. The computations   were carried out using MOSEK solver \citep{mosek} and YALMIP interface \citep{Lofberg2004} on MATLAB.
\subsection{Bounds for various bivariate dependencies}
In the first set of experiments, for $n=15$ Bernoulli random variables, we considered $50$ randomly generated trees with the univariate probabilities $p_i$  generated uniformly in $(0,0.1]$. We computed the following bounds for various values of $k$.
(1) $U_{uv}(k) $ = Maximum value of $\mathbb{P}(\sum_{i=1}^n \tilde{c}_i \geq k)$ assuming univariate information alone (using \eqref{eqn:univar-bound})
(2) $U(k)$ = Maximum value of $\mathbb{P}(\sum_{i=1}^n \tilde{c}_i \geq k)$ assuming bivariate distributional information (using \Cref{bivar-tree-theorem}) (3)$P_{ci} = \mathbb{P}(\sum_{i=1}^n \tilde{c}_i \geq k)$ for the conditionally independent distribution (using \Cref{thm:cond_ind_recursion}) for the same bivariate information used in the computation of $U(k)$.  We study the scenarios when the bivariate probabilities are generated using the comonotone and anti-comonotone copulas with the generated univariate probability distributions. The comonotone copula represents maximum  positively dependent random variables while the anti-comonotone  copula represents maximum negatively dependent random variables \citep{Nelsen2006,puccetti2015}.  For Bernoulli random variables, it is known that for the comonotone copula, $p_{ij}= \min (p_i, p_j)$ while for the anti-comonotone copula, $p_{ij}= \max(p_i + p_j - 1, 0)$.

\begin{figure}[h]
\centering
\caption{Demonstration of the bounds for various dependence structures. The  range of bounds $U, U_{uv}$ and $P_{ci}$ over $50$ arbitrary tree structures on $n=15$ random variables is shown.}
\label{fig:bivar-univar-bound}
\begin{subfigure}{0.45\textwidth}
\includegraphics[scale=0.55]{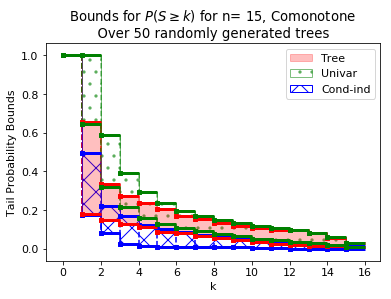}
\end{subfigure}
\begin{subfigure}{0.45\textwidth}
\includegraphics[scale=0.55]{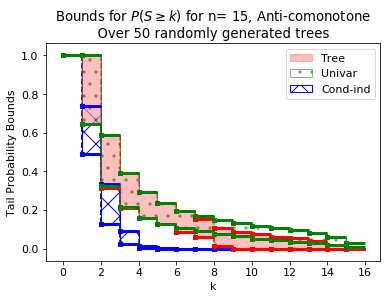}
\end{subfigure}
\end{figure}

The range of values of $U_{uv}$, $U$ and $P_{ci}$ over $50$ runs are provided in Figure \ref{fig:bivar-univar-bound} (with labels Univar, Tree and Cond-ind respectively). The dotted, shaded and cross-hatched region correspond to the univariate, tree and conditionally independent bounds respectively. We see that when the bivariate distributions are specified using the comonotone copula, the univariate bound is much larger than the tree bound for lower values of $k$ whereas for larger values of $k$ they almost co-incide. On the other hand, for the anti-comonotone copula, the tree bound is almost identical to the univariate bound for lower values of $k$ and the values start differing as $k$ becomes large. This makes intuitive sense as when no constraints on bivariates are specified, we expect the largest value of $\mathbb{P}(\sum_{i=1}^n \tilde{c}_i \geq k)$ to be attained by the comonotone copula for larger values of $k$ as the comonotone copula assigns larger probabilities for more random variables taking similar values. Thus the distribution that achieves  the optimal univariate bound drifts closer to the comonotone distribution and anti-comonotone distributions for larger and smaller values of $k$ respectively.
Note that in general, the conditionally independent distribution gives a different bound. As it is a feasible distribution in $\Theta$, it gives a lower bound to the optimal tree and univariate bounds.

\subsection{Robustness of the Conditionally Independent Chow-Liu Tree}
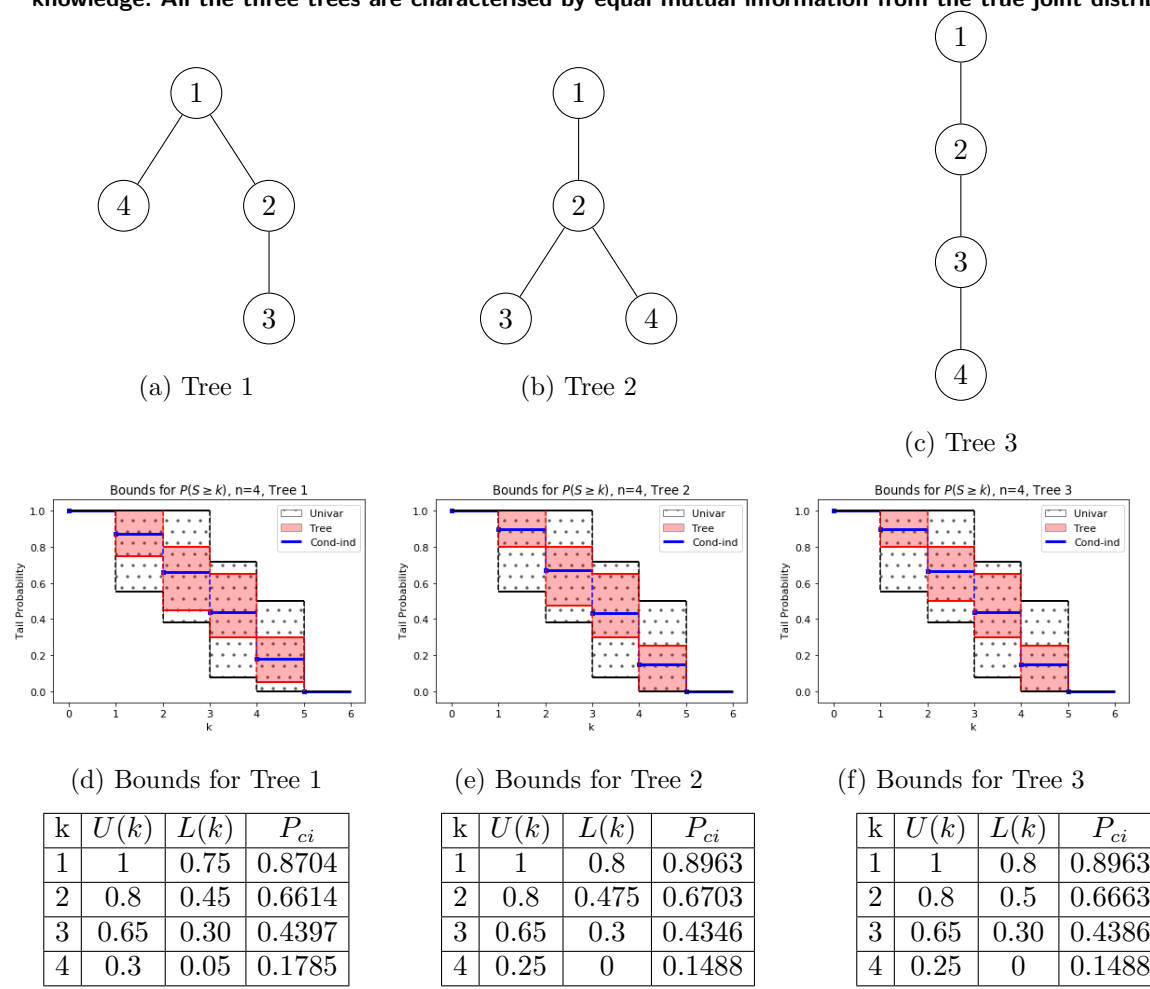
\begin{figure}[h!]
\caption{The trees in    \cite{chowliutree1968} and the bounds produced on $\mathbb{P}(\tilde{c}_i \geq k)$ under various assumptions of knowledge. All the three trees are characterised by equal mutual information from the true joint distribution. }
\label{fig:CLTree-bounds}
\begin{subfigure}{0.3\textwidth}
\centering
\begin{tikzpicture}[sibling distance=5em,
  every node/.style = {shape=circle,
    draw, align=center,
    top color=white}]]
  \node {1}
    child { node {4} }
    child { node {2}
      child { node {3} }};
\end{tikzpicture}
\caption{Tree 1}
\end{subfigure}
\begin{subfigure}{0.3\textwidth}
\centering
\begin{tikzpicture}[sibling distance=5em,
  every node/.style = {shape=circle,
    draw, align=center,
    top color=white}]]
  \node {1}
    child { node {2}
    child { node {3}}
      child { node {4} }};
\end{tikzpicture}
\caption{Tree 2}
\end{subfigure}
\begin{subfigure}{0.3\textwidth}
\centering
\begin{tikzpicture}[
  every node/.style = {shape=circle,
    draw, align=center,
    top color=white}]]
  \node {1}
    child { node {2}
    child { node {3}
      child { node {4} }}};
\end{tikzpicture}
\caption{Tree 3}
\end{subfigure} \\
\begin{subfigure}{0.3\textwidth}
\includegraphics[scale=0.35]{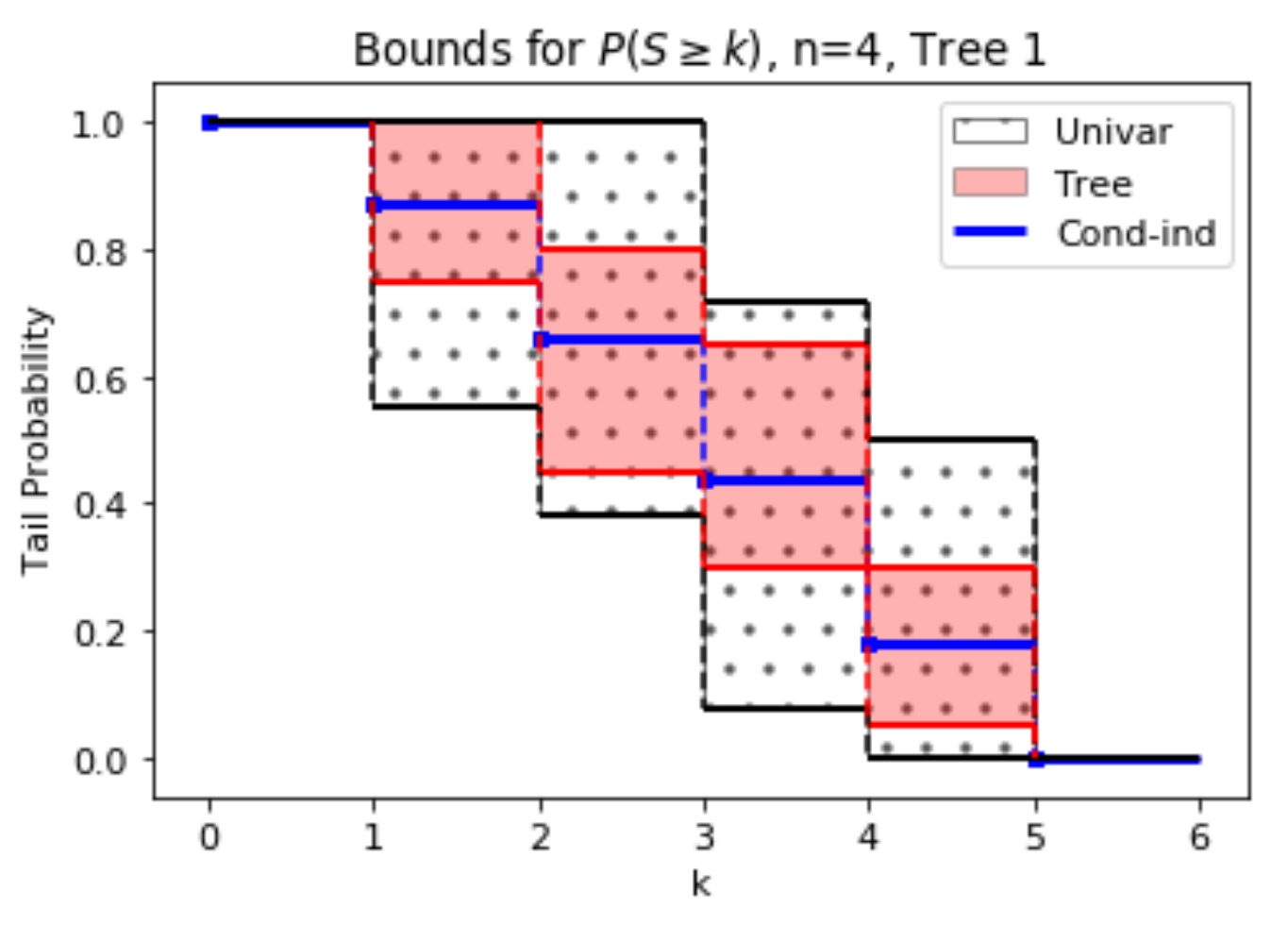}
\caption{Bounds for Tree 1}
\end{subfigure}
\begin{subfigure}{0.3\textwidth}
\includegraphics[scale=0.35]{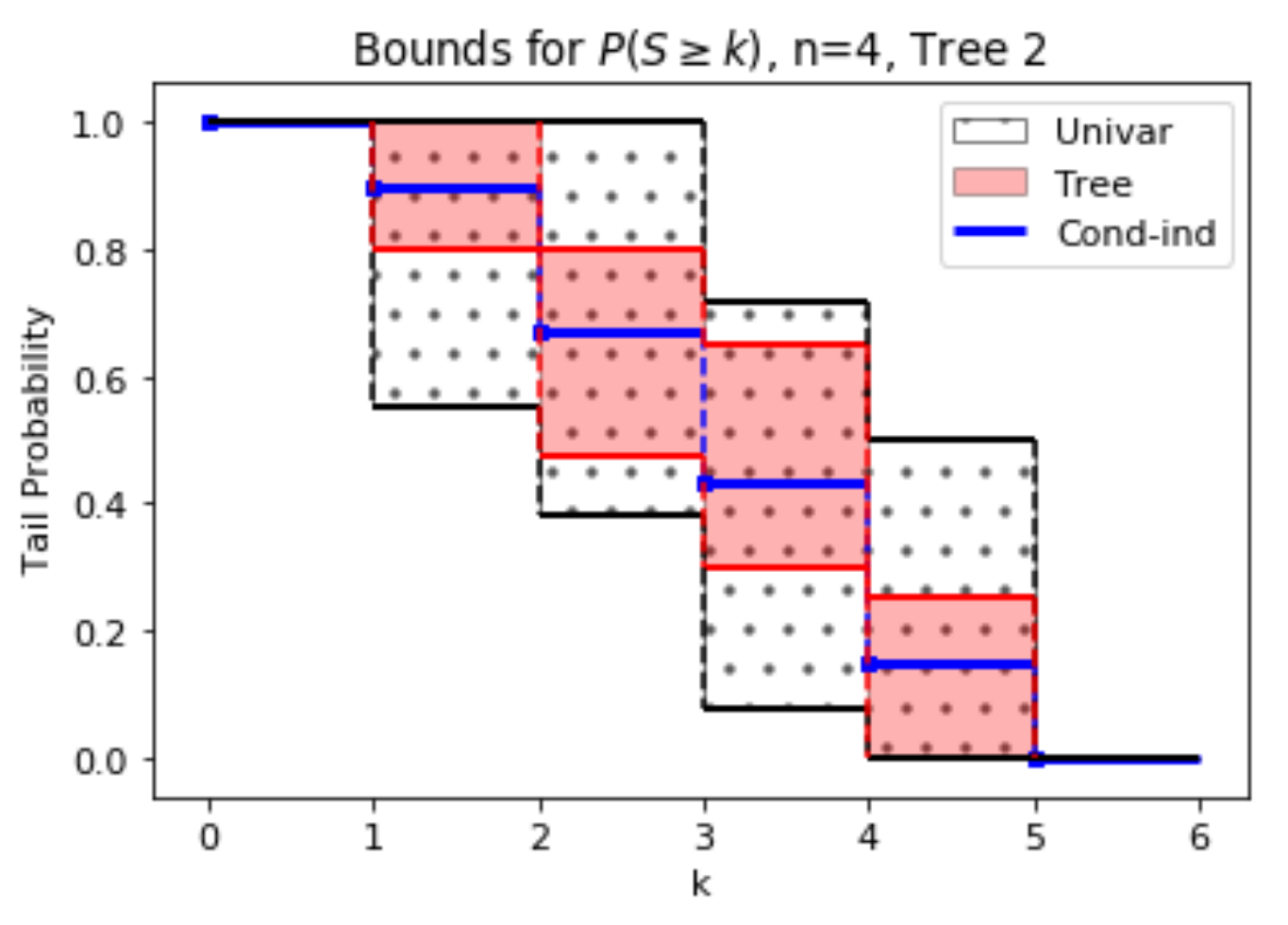}
\caption{Bounds for Tree 2}
\end{subfigure}
\begin{subfigure}{0.3\textwidth}
\includegraphics[scale=0.35]{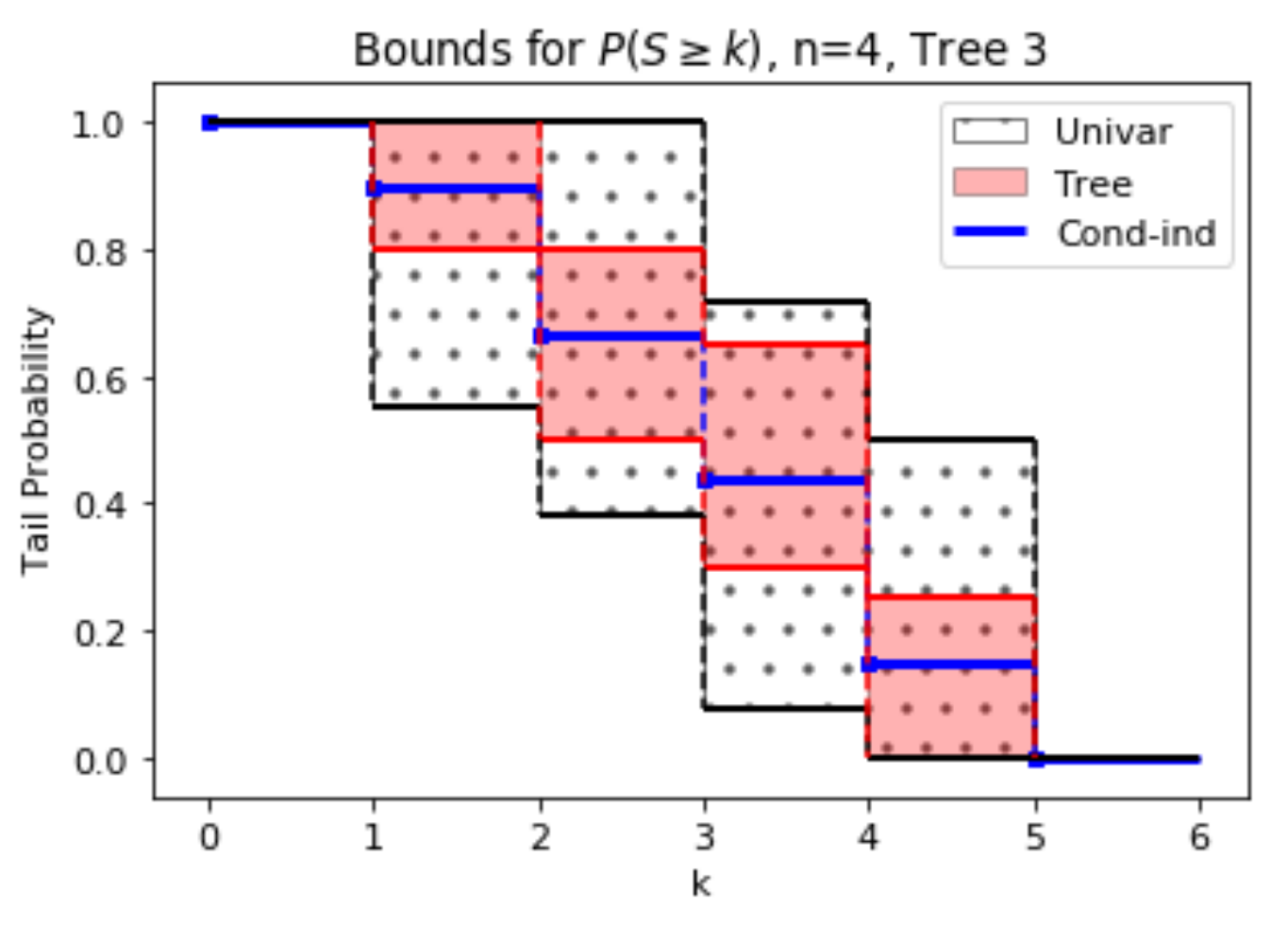}
\caption{Bounds for Tree 3}
\end{subfigure}\\
~~~\begin{minipage}{0.3\textwidth}
\centering
\begin{tabular}{|c|c|c|c|}
\hline
k & $U(k)$ &$L(k)$ & $P_{ci}$\\\hline
1	& 1  &   0.75 &  0.8704   \\\hline
2	& 0.8  &  0.45 & 0.6614  \\\hline
3	& 0.65  &   0.30 & 0.4397   \\\hline
4	& 0.3  &  0.05  &  0.1785  \\\hline
\end{tabular}

\end{minipage} ~
\begin{minipage}{0.3\textwidth}
\centering
\begin{tabular}{|c|c|c|c|}
\hline
k & $U(k)$ &$L(k)$ & $P_{ci}$\\\hline
1	& 1  &   0.8  &  0.8963 \\\hline
2	& 0.8  &  0.475  & 0.6703 \\\hline
3	& 0.65  &   0.3   & 0.4346 \\\hline
4	& 0.25  &  0   & 0.1488 \\\hline
\end{tabular}
\end{minipage} ~~
\begin{minipage}{0.3\textwidth}
\centering
\begin{tabular}{|c|c|c|c|}
\hline
k & $U(k)$ &$L(k)$ & $P_{ci}$\\\hline
1	& 1  &   0.8 & 0.8963   \\\hline
2	& 0.8  &  0.5 & 0.6663 \\\hline
3	& 0.65  &   0.30 & 0.4386 \\\hline
4	& 0.25  & 0 & 0.1488 \\\hline
\end{tabular}
\end{minipage}\\\\
\end{figure}

In this subsection we consider the probability distribution on $n=4$ Bernoulli random variables provided in  \cite{chowliutree1968}, Table 1.  The univariate probabilities  computed from the joint distribution in the original paper are $p_1 = p_2 = p_3 = 0.55, p_4 = 0.5$ while the bivariate probabilities are $p_{12} = 0.4, p_{13} = p_{14} = 0.3, p_{23} = 0.45, p_{24} = 0.25,  p_{34} = 0.25$. Three trees that best approximate the provided joint distribution were also provided in the same paper. We present the trees themselves in \Cref{fig:CLTree-bounds}. The trees are equivalent in that the sum of mutual information encoded by the probability distributions on the edges of the trees is the same.

Given a tree, the most natural distribution is the conditionally independent joint  distribution. As the tree structures are different, the conditionally independent distributions themselves differ for these three trees.  We report the following set of bounds:
\begin{enumerate}
\item The probability $P_{ci}= \mathbb{P}(\sum_{i=1}^n \tilde{c}_i \geq k)$ assuming a   conditional independent joint distribution on each of these trees (using \Cref{thm:cond_ind_recursion}).
\item  Upper and lower bounds  $U(k)$ and $L(k)$ on $\mathbb{P}(\sum_{i=1}^n \tilde{c}_i \geq k)$ assuming the bivariate distributions $\mathbb{P}(\tilde{c}_i = 1,\tilde{c}_j = 1 )$ for all edges $(i,j)$ in the given tree (using \Cref{bivar-tree-theorem} and \Cref{eq:lower_bound} respectively).
\item Upper and lower bounds  $U_{uv}(k)$ and $L_{uv}(k)$ on $\mathbb{P}(\sum_{i=1}^n \tilde{c}_i \geq k)$ assuming the univariate probabilities $\mathbb{P}(\tilde{c}_i = 1)$ alone. $U_{uv}(k)$ was computed using \Cref{eqn:univar-bound} while, by a similar argument in \Cref{eq:lower_bound}, $L_{uv}(k) = U_{uv}(n-k+1)$ (where, in the computation of $U_{uv}(n-k+1)$ the probabilities $1-p_i$ were sorted and used instead of $p_i$).
\end{enumerate}

The univariate probabilities thus  computed are provided in \Cref{tab:chow-liu-univar}. These probabilities are the same for all the trees as they do not make use of any bivariate information. The range $[L_{uv}(k), U_{uv}(k)]$ is provided in the plots in \Cref{fig:CLTree-bounds} as the dotted region. This band  forms the widest band as the bounds assume only univariate information. Under assumptions of bivariate information on the trees, we see that the band $[L(k), U(k)]$ (plots depicted by the shaded region in  \Cref{fig:CLTree-bounds}) is narrower as more information is assumed. These bands are clearly contained in the univariate band. In each of the plots, we see that the probability computed using the conditionally independent distribution  (shown in blue) is much farther away from the upper as well as lower bounds and sits in the middle of the bands. This demonstrates  examples where the conditional independence on the Chow-Liu tree may approximate the provided distribution well, however the optimal upper and lower bounds on the tail probabilities are achieved by a different distribution.

\begin{table}[h]
\centering
\begin{tabular}{|c|c|c|}
\hline
k & $U_{uv}(k)$ &$L_{uv}(k)$ \\\hline
1	& 1  &   0.55    \\\hline
2	& 1  &  0.3833 \\\hline
3	& 0.7167  &   0.30 \\\hline
4	& 0.5  & 0  \\\hline
\end{tabular}
\caption{Probability bounds assuming only univariate information}
\label{tab:chow-liu-univar}
\end{table}

\subsection{Bounds for Order Statistics Probabilities}
In this subsection we present various bounds for  the order statistics probabilities. We consider the following two distributions

(1) A multi-variate Gaussian distribution in $5$ dimensions with  randomly generated mean $\mu$  provided in \Cref{gaussian-mean} and covariance matrix $\Sigma $ taken as the identity matrix $\mathbbm{I}_5$, in 5 dimensions.
\begin{align}
\mu & =\begin{bmatrix}&  0.5426, & -0.9585, &
    0.2673, &
    0.4976, &
   -0.0030\end{bmatrix}   \label{gaussian-mean}
\end{align}

(2) A multi-variate Pareto distribution  in $5$ dimensions  with the following parameters, shape = $1$ and location $l_{Pareto} = [6.54, 5.04, 6.26, 6.49, 5.99]$.

\begin{figure}[h]
\begin{subfigure}{0.33\textwidth}
\includegraphics[scale=0.4]{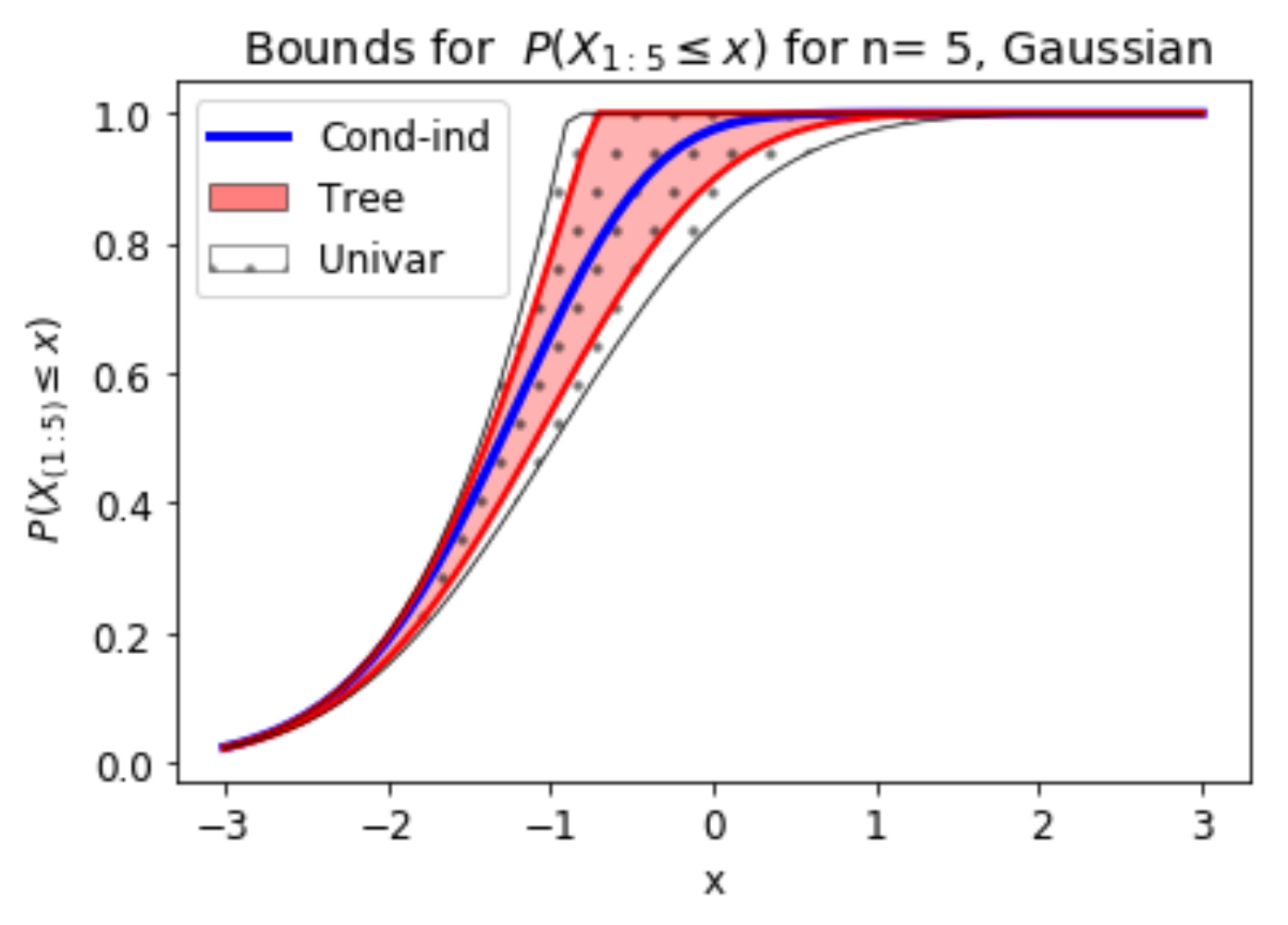}
\subcaption{Min Order Statistic\\  (Gaussian)}
\label{fig:gaussian_min_orderstat}
\end{subfigure}
\begin{subfigure}{0.33\textwidth}
\includegraphics[scale=0.4]{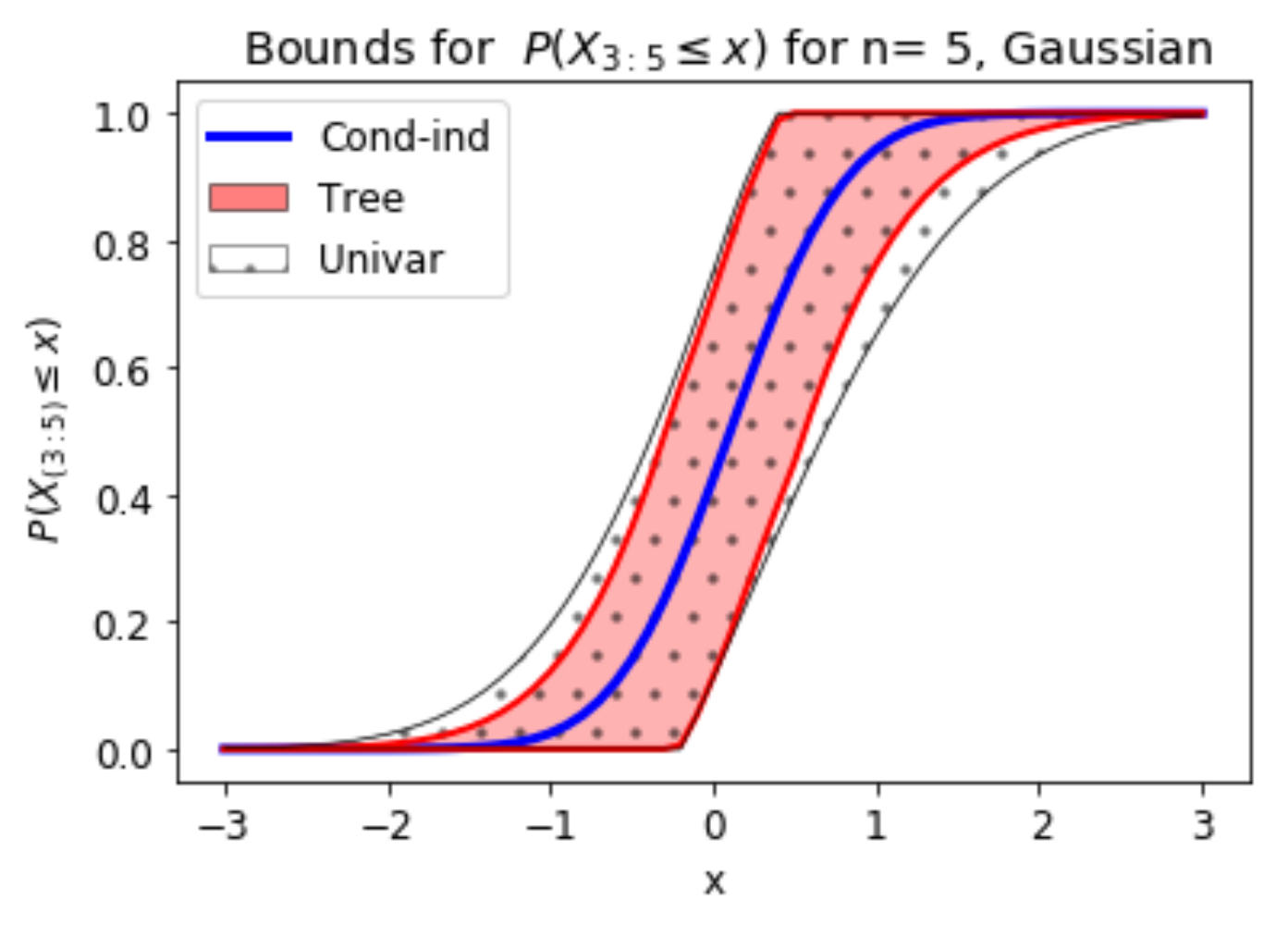}
\subcaption{Median Order Statistic\\  (Gaussian)}
\label{fig:gaussian_median_orderstat}
\end{subfigure}
\begin{subfigure}{0.33\textwidth}
\includegraphics[scale=0.4]{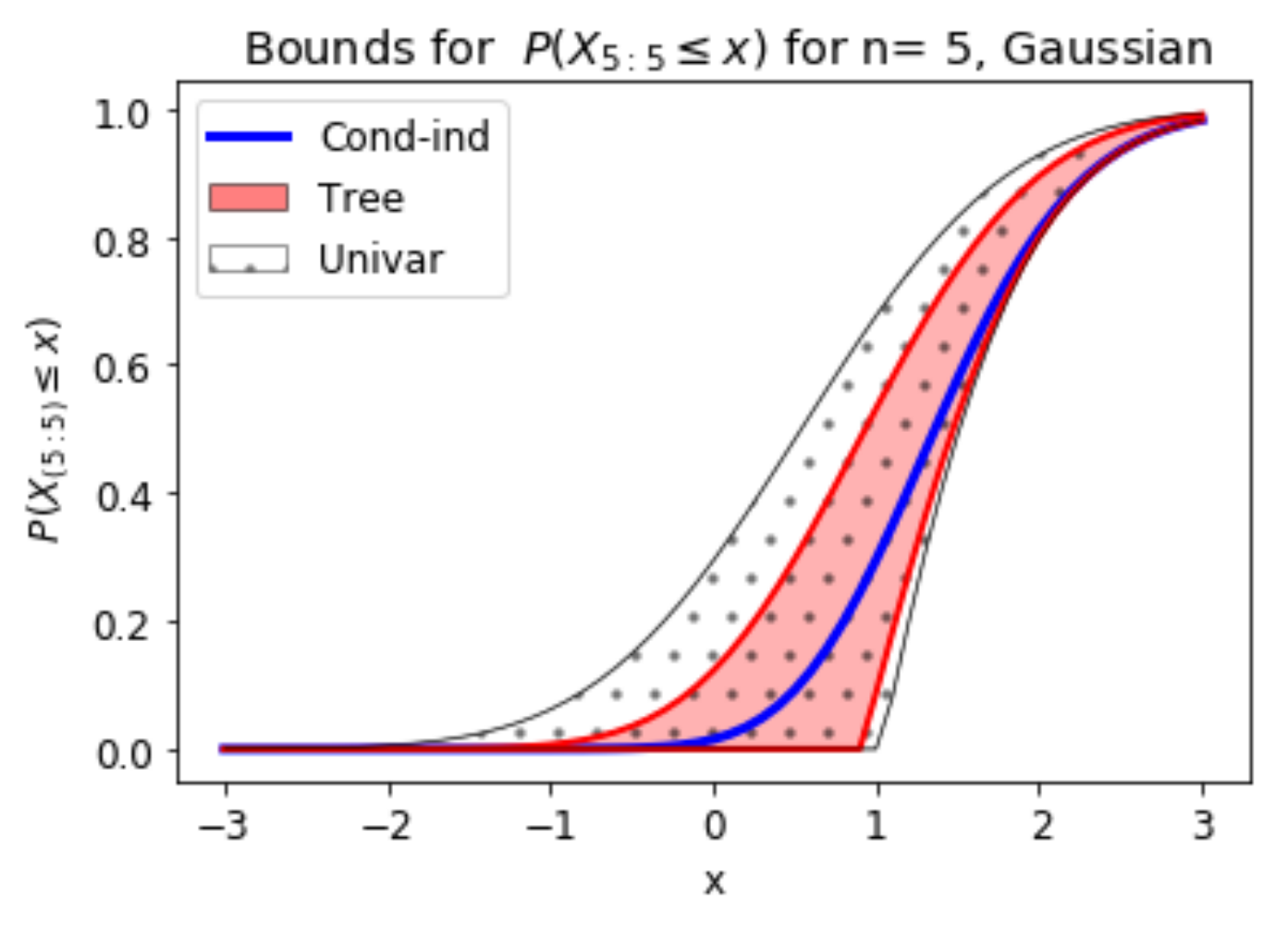}
\subcaption{Max Order Statistic\\  (Gaussian)}
\label{fig:gaussian_max_orderstat}
\end{subfigure} \\
\begin{subfigure}{0.33\textwidth}
\includegraphics[scale=0.4]{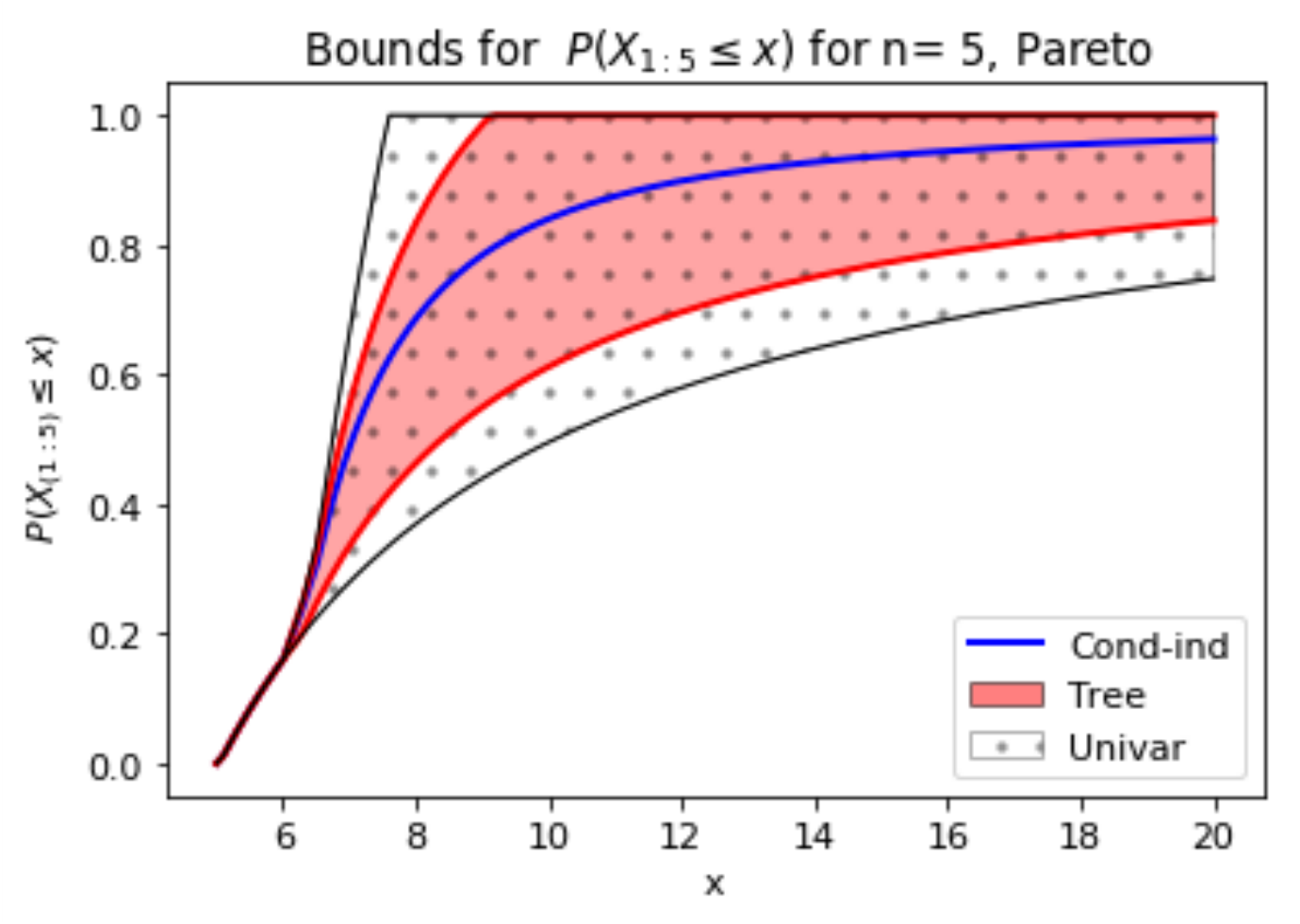}
\subcaption{Min Order Statistic (Pareto)}
\label{fig:pareto_min_orderstat}
\end{subfigure}
\begin{subfigure}{0.33\textwidth}
\includegraphics[scale=0.4]{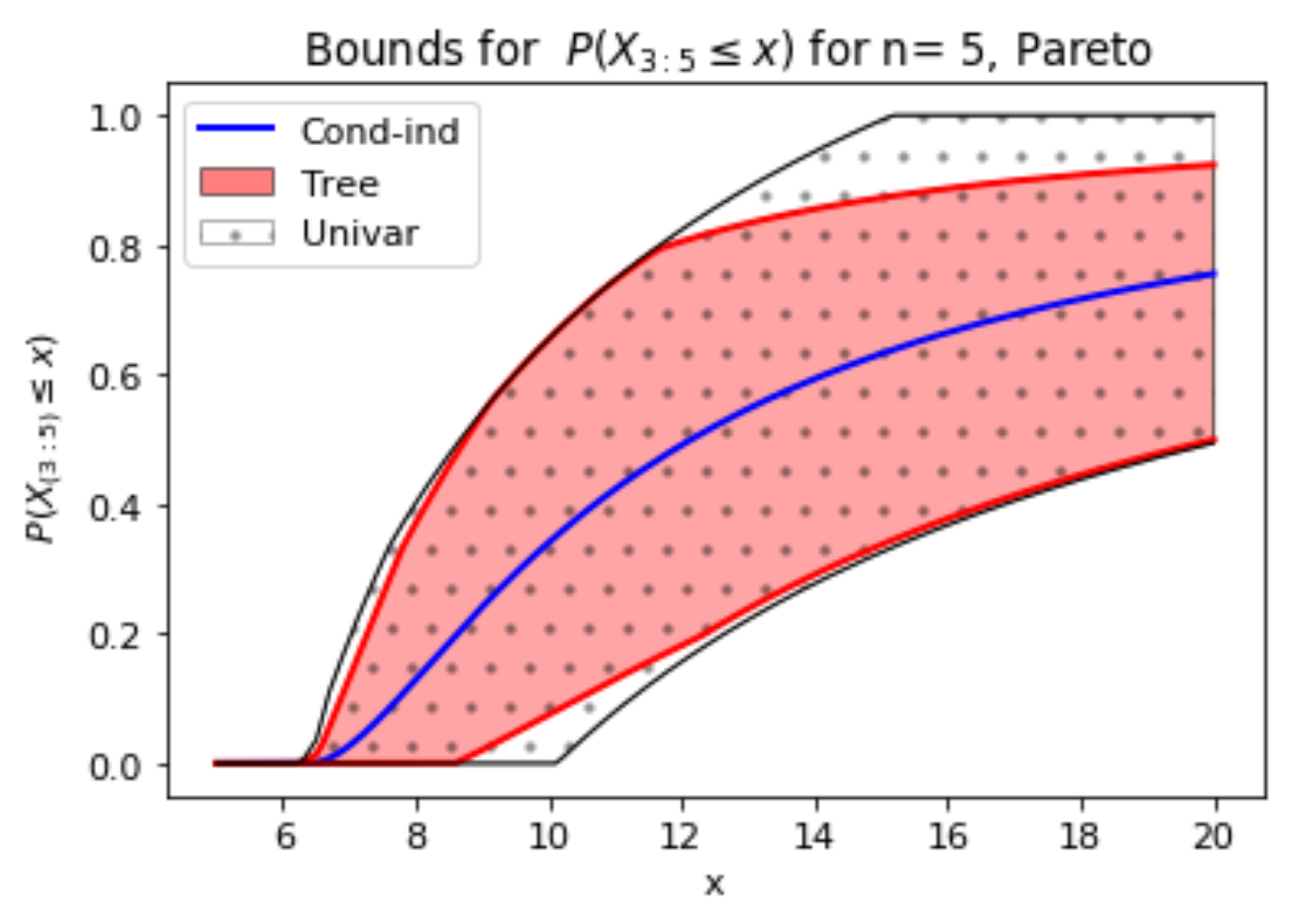}
\subcaption{Median Order Statistic (Pareto)}
\label{fig:pareto_median_orderstat}
\end{subfigure}
\begin{subfigure}{0.33\textwidth}
\includegraphics[scale=0.4]{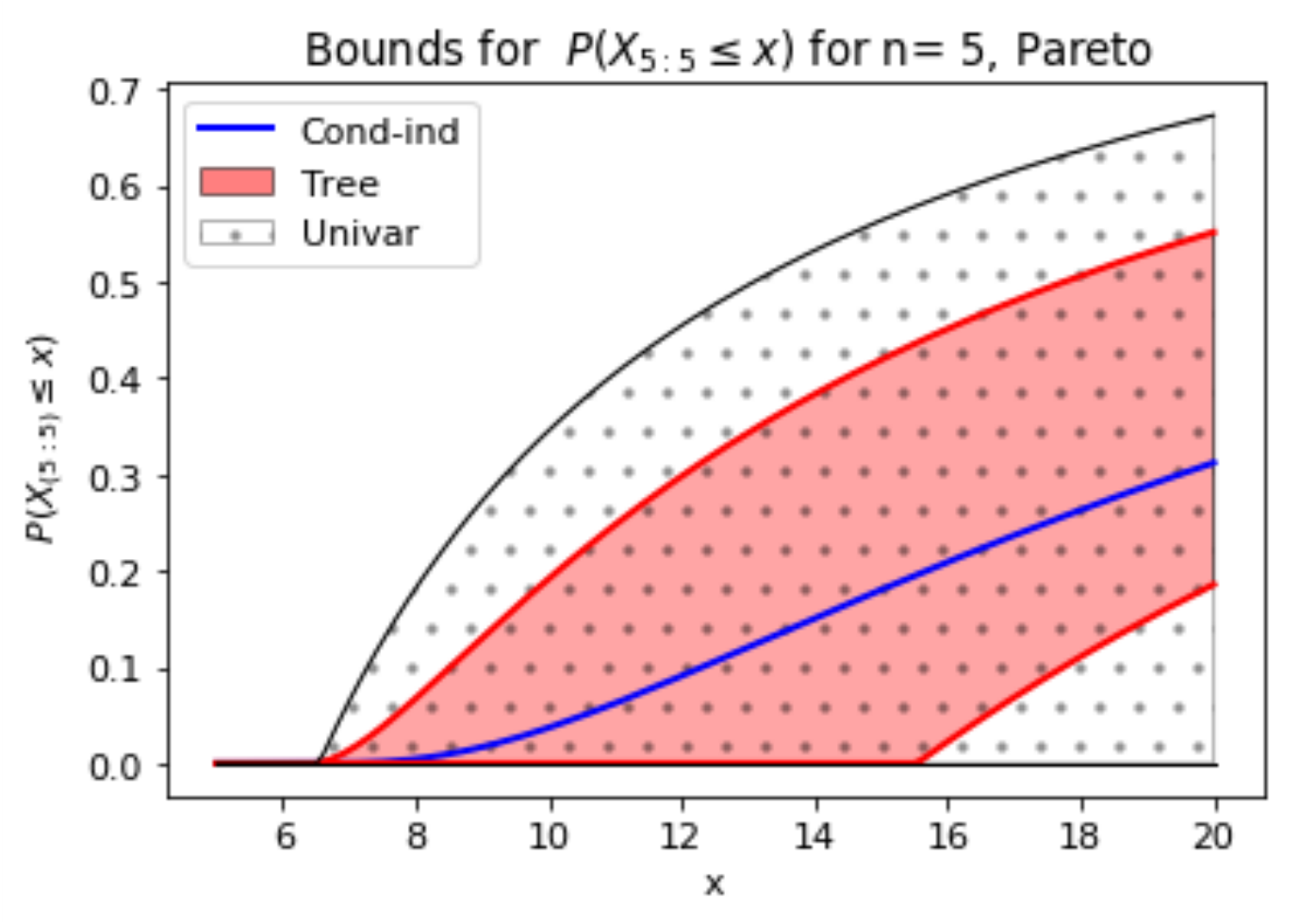}
\subcaption{Max Order Statistic(Pareto)}
\label{fig:pareto_max_orderstat}
\end{subfigure}
\caption{Bounds for the $k^{th}$ order statistic probabilities}
\label{fig:order-stats-prob}
\end{figure}

The Pareto distribution lies in the class of heavy tailed distributions while the multi-variate Gaussian is an example of a light tailed distribution. The bounds on the $k^{th}$ order statistics probabilities were computed for $x \in [-3, 3]$ in steps of $0.1$ and $k \in \{1, \ldots, 5 \}$ for the Gaussian distribution, while for the Pareto distribution the range $x \in [5, 20]$ was used. We illustrate the following in  \Cref{fig:order-stats-prob}:
\begin{enumerate}
\item The range $[L_{os,k}(x),U_{os,k}(x)]$ for the order statistic probabilities using  \Cref{thm:order-stats-bd} and the relations between the upper and lower bounds in \Cref{sec:lower_bound}, assuming bivariate distributional information $\mathbb{P}(\tilde{X}_i \leq x, \tilde{X}_j \leq x)$ on the series graph: For each order statistic $X_{k:n}$, the minimum and maximum value of $\mathbb{P}(X_{k:n} \leq x)$ for various values of $x$ were computed and the range of $\mathbb{P}(X_{k:n} \leq x)$ is represented as the shaded region (labelled `Tree').

\item  Order statistic probabilities using  \Cref{eqn:univar-bound}, assuming only univariate distributional information and  setting $p_i  = \mathbb{P}(\tilde{X}_i \leq x)$. For a given $x$, $L_{os,k}^{uv}(x)$ and  $U_{os,k}^{uv}(x)$ denote the minimum and maximum values of $\mathbb{P} (X_{k:n} \geq x)$ computed under assumptions of univariate distributional information alone. Analogous to the case of bivariate information, $L_{os,k}^{uv}(x) = U_{os,n-k+1}^{uv}(x)$ where in the computation of $U_{os,n-k+1}^{uv}(x)$, all operations involving $p_i$ are replaced with $1-p_i$. For various values of $x$, the range $[L_{os,k}^{uv}(x), U_{os,k}^{uv}(x)]$  is  provided by the dotted region (labelled `Univar').

\item Order statistic probabilities for the conditionally independent distribution on the series graph using \Cref{thm:cond_ind_recursion}. This cumulative density function is shown by the blue line labelled `Cond-ind'.

\end{enumerate}

The min, median and max order statistic probabilities for the Gaussian and Pareto distributions are  shown in \Cref{fig:gaussian_min_orderstat,fig:gaussian_median_orderstat,fig:gaussian_max_orderstat}
and \Cref{fig:pareto_min_orderstat,fig:pareto_median_orderstat,fig:pareto_max_orderstat}
respectively.  As expected in all cases, the  tree band  is sandwiched in the univariate band. This is natural as the tree bounds make use of more information than the univariate bounds. The probabilities given by the conditionally independent distribution  lie in the region spanned by the tree bounds. In general, we observe that the lower and upper tree bounds are very different from the bound provided by the conditionally independent distribution indicating that in scenarios where robustness or extremal values are of interest, the conditionally independent distribution is not the best.

Notice that the minimum and maximum value of the support for each of the bounds considered progressively becomes larger as we go from the min to the max order statistic. For example, for the Gaussian distribution, the support for the upper tree bound ($\approx [-3,-1]$) for the min order statistic moves to  $\approx [-2,0]$ for the median order statics  and then to $\approx  [-1, 3]$ for the maximum order statistic. Similar trend is exhibited by all other bounds. Also the the range for this support is wider for the case of Pareto distribution. For example, the support for the upper tree bound for the min order statistic range $\approx$ [5,9] and  [6, y] for the median and max order statistic of the Pareto distribution,  where $y \geq 20$. This is a consequence of the heavy tail behaviour of the Pareto distribution which is well exhibited by our tree bounds too.

\clearpage
\small
\bibliographystyle{apalike}
\bibliography{tree-bivar}

\end{document}